\documentclass[12pt]{article}	

\usepackage[utf8]{inputenc}				
\usepackage{epsfig}
\usepackage{amssymb,amsmath,amsfonts,amsthm,enumerate}
\usepackage{amsfonts}
\usepackage{amsmath}
\usepackage{graphicx}
\usepackage{amsthm}
\usepackage{mathrsfs}
\usepackage{color}

\begin{document}


\begin{center}
    \bfseries
		
{\large Some aspects of neural network parameter optimization for joint inversion of gravitational and magnetic fields}

    \bigskip
    Yanfei Wang\footnotemark[1], Dmitry V. Churbanov\footnotemark[2], Raul L. Argun\footnotemark[3], Alexander V. Gorbachev\footnotemark[3],  Alexander S. Leonov\footnotemark[4], Dmitry V. Lukyanenko\footnotemark[3]

\footnotetext[1]{Institute of Geology and Geophysics, Chinese Academy of Sciences, Beijing 100029, China, e-mail: yfwang@mail.iggcas.ac.cn}
    \footnotetext[2]{Moscow Institute of Physics and Technology, Department of Higher Mathematics, Moscow 141700, Russia, e-amil: churbanov.dv@phystech.edu}
    \footnotetext[3]{Lomonosov Moscow State University, Faculty of Physics, Department of Mathematics, Moscow 119991, Russia, e-mail: argun.rl14@physics.msu.ru,  gorbachev.av17@physics.msu.ru, lukyanenko@physics.msu.ru}

    \footnotetext[4]{National Research Nuclear University MEPhI, Department of Higher Mathematics, Moscow 115409, Russia, e-amil: asleonov@mephi.ru}

\end{center}


\bigskip	
{\small\textbf{Abstract}. We consider the optimization of a neural network previously developed by the authors for the joint inversion of 3D gravitational and magnetic fields in the context of mineral exploration. The distinctive feature of this neural network is that it solves ill-posed (ill-conditioned) inverse problems. The neural network implements a special two-level algorithm. The lower level of the algorithm uses two neural networks with equivalent architectures. The first of them computes the gravitational field sources in a given domain from measurements of this field on a remote surface. The second neural network processes magnetic field measured on the same surface to find magnetic sources in the same domain. The found source distributions are used at the upper level of the algorithm to calculate their structural residual, which determines the degree of difference (closeness) of their geometries. As a result, minimizing this residual, when training a neural network at the upper level, implements a computational algorithm that yields geometrically close source distributions of different fields. The article examines in detail the possibilities of optimizing some elements of the neural networks and the algorithms used (datasets, training process, specific form of loss functions, etc.) Test calculations for model problem demonstrate high quality of joint inversion by our optimized neural networks approach. Calculations were also carried out for the joint processing of real-feald data from gravity and magnetic exploration in Jussara region, Goias State, Brazil. The article also considers the issue of determining in joint field inversion not only the geometric distribution of sources, but also their physical intensities.}
\bigskip

{\bfseries Keywords:} inverse problem, joint inversion of gravitational and magnetic fields, neural networks, structural residual.

\section{Introduction}

In recent years, there has been significant interest in the application of neural networks for solving various inverse problems. Neural networks of different architectures, trained by various methods, are increasingly used to solve typical inverse problems of varying complexity in science, engineering, economics, etc. This is largely explained by the fact that typical inverse problems often have strictly characteristic solutions. Therefore, the following, not very detailed, scheme for the approximate solution of the corresponding inverse problem using a neural network of a certain architecture becomes practically acceptable.
\begin{enumerate}
 	\item Creation of a sufficiently representative ensemble of characteristic solutions for training (dataset), i.e., a set of pairs (data, solutions).
	\item Specification of the target function (loss) for network training, i.e., a mathematical criterion for the proximity of the data presented to the network and their analogs computed by the network.
	\item Training a neural network with a given architecture and loss function on the selected dataset, i.e., from a mathematical point of view, constructing an operator for solving the inverse problem for data from the dataset.
	\item Validation of the trained neural network.
	\item Processing new data of the inverse problem (usually not included in the dataset) using the trained neural network, i.e., applying the found inverse problem solution operator to this new data.
\end{enumerate}

Such a scheme allows, with an adequate choice of network architecture, dataset, loss criterion, and training procedure, to quickly and relatively accurate solve the inverse problem, at least for data close to the data from the dataset.

An obvious drawback of this scheme is its possible inaccuracy for data not belonging to the dataset, as well as the instability of the obtained solutions with respect to small variations in the data that arises when the neural network is overfitted. Accordingly, a ready-to-use neural network must be properly ``optimized''. We understand this as meaning that for a characteristic class of model inverse problems, it should yield solutions that are stable to data perturbations with accuracy acceptable to the user. Such optimization will be investigated in our article for special inverse problems of joint inversion of gravitational and magnetic fields. The optimization will be associated with taking into account \textit{a priori} information about the fields and solutions in the problems of separate inversion. It is carried out by special construction of the loss function, the dataset, and by introducing additional solution quality functionals into the inverse problem.

Let us say a few general words about the inverse problems considered below. Problems of potential field inversion often arise in gravity and magnetic prospecting, being an important part of the practical search for minerals. Among such problems, those related to the restoration of the shape and location of three-dimensional ore structures are of significant interest. The data here are various measurements of gravitational and magnetic fields, their gradients, etc. The search for distributions of gravitational (magnetic) field sources from this data is commonly called the inverse problem of gravimetry (magnetometry). Historically, the corresponding mathematical inverse problems were posed and solved separately~\cite{L1,L2,TA,Zhdanov-1984}. However, recently there has been great interest in their joint solution (see, for example, the review~\cite{simirdanis_report}).

Each of these inverse problems is ill-posed, and when solving them separately, Tikhonov regularization is often used~\cite{TA,TLY,L5}. The work~\cite{Wang-2019-2} can serve as an example. Methods are also used in which regularization is provided by sufficiently detailed \textit{a priori} information about the solution. \textit{A priori} information can also be taken into account to simplify the solution procedures. For example, if it is known that the reconstructed subsurface object has a layered structure with a known density for each layer, then one can pose the inverse problem of reconstructing the shape of the surfaces separating these layers (see, for example,~\cite{Akimova-2023, huang22}). In this case, the original problem of restoring the three-dimensional domain of the ore body is reduced to a two-dimensional problem of determining the interface surfaces. The resulting reduction in the dimensionality of the problem significantly reduces its computational complexity. If there is no significant \textit{a priori} information about the solution, then the problems of potential field inversion have to be solved in the general (three-dimensional) formulation. This, in turn, is associated with significant theoretical and computational difficulties. Even greater difficulties arise when attempting {joint field inversion}.

In this work, we consider the problem of joint inversion of potential fields, assuming that the sources of the gravitational and magnetic fields are concentrated in the same region and have the same parameterized geometric shape. Our goal is to find the spatial structure and location of a \emph{single ore body} consisting of these sources. The data for solving the inverse problem in this work are measurements of potential fields (gravitational and magnetic) created by the ore body on the terrain surface or near it. Our proposed approach to solving this problem will be based on the use of neural networks and machine learning methods~\cite{Yu21}. This approach opens up opportunities for the development of effective methods specifically in the \emph{three-dimensional joint inversion} of potential fields. The basis for this is the expressive ability of neural networks, i.e., the ability to effectively reproduce complex dependencies.

As noted above, when designing and using neural networks, in particular for the problems we consider, it is necessary to solve several problems: choose the type and architecture of the neural network, choose training parameters, set the training criterion (loss function), and adequately define the set of training data (dataset). These choices should reflect the specifics of the problem being solved, and this will be discussed below. But for the problem of joint inversion of potential fields using neural networks, a number of additional problems arise. Most of them are related to the ill-posedness of the inverse problem being solved.

First of all, this is the non-uniqueness of its solution (see~\cite{L1,L2,L3}). The same set of measurements of the gravitational and/or magnetic field on a given surface can correspond to different sets of three-dimensional subsurface sources in terms of geometric configurations and physical properties. Another problem is the instability of the solutions of the inversion problem with respect to small perturbations of the data. ``Close'' data (in the considered metric) can correspond to significantly different subsurface structures. This can lead to unacceptable errors in determining such characteristic features of ore structures as faults and fractures. Finally, the third problem is the impossibility of assessing the accuracy of the obtained solutions of the inverse problem without having detailed \textit{a priori} information about the solution. The problem is related precisely to the ill-posedness of the inverse problem being solved and it is well known from theory (see, for example,~\cite{L4,L5}).

In principle, all these problems are solved in conventional (non-neural network) approaches to field inversion problems by finite-dimensional parameterization of their sources with known density or magnetic properties of the substance. Then the sources form a compact set in the corresponding functional space, and the solution of the inverse problem often turns out to be unique and stable with respect to data perturbations. On a compact class of solutions, one can also obtain an estimate of their accuracy (see, for example,~\cite{IVT,L4,L5}). This idea is convenient to apply also in the neural network approach through the use of a specialized dataset.

When training neural networks, not only the sufficient volume of the dataset is very important, but also its qualitative composition (in terms of the types of representatives). It should include data and results of joint inversion specifically for the used parametric classes of sources. It is also desirable to include the results of processing geophysical data there, although the number of such results in the public domain is very small. Accordingly, to build an adequate dataset, it is necessary to partially generate solutions of inverse problems from models and, if possible, supplement such artificial data with real ones. In the general case, with having minimal \textit{a priori} information about the solution in joint inversion of fields created by bodies of various irregular shapes, one has to use only a \emph{synthetic dataset}. Its parameterization significantly affects the stability of the solutions of the inverse problem. It is convenient to include such parameterized forms of potential field sources as combinations of prisms and steps of various shapes with known substance density, as well as bodies of irregular shape generated based on some stochastic process from a random set of smaller cubes (see, for example,~\cite{huang21}). Other ways of obtaining a dataset from various three-dimensional bodies of irregular shape are also possible (see~\cite{qiao24}). In any case, such a dataset must take into account the available \textit{a priori} information about the solution. In particular, this can be statistical information from expert geologists, which includes various geological scenarios and estimates (see~\cite{noddyverse21}).

Various attempts are being made to optimize the dataset, which allow improving the quality of the solution while preserving the original model. Here, one of the most universal and efficient approaches in terms of resources spent is augmentation~\cite{zhang21dive}. In this approach, the used dataset is subjected to various modifications, which increase its diversity without the costly generation of new objects. For example, the work~\cite{zhou23} implements dataset augmentation for gravitational inversion. The initial sample there consists of simple geometric figures: prisms and steps. They are subjected to reflection transformations, addition of original forms with a proportionality coefficient, and noising. This essentially uses the linearity of the problem, so that data for new elements created by augmentation in the dataset are easily computed without solving the full forward problem. Accordingly, an analysis of the influence of various types of augmentation on the final solution accuracy is carried out. Another possible approach to dataset optimization is relevant for problems of gravitational or magnetic inversion when processing data for ore bodies of complex shape and large sizes. For example, the work~\cite{huang22} proposes to break large problems into several small ones, using the linearity property of the original forward gravimetry problem. The sought ``large'' body of complex shape is divided into several smaller subdomains of simpler shape, which are easier to reproduce by the neural network given current datasets with small geometric sizes of solutions. Therefore, building datasets of large geometric sizes, which is usually difficult from a computational point of view, is not required.

The training of any neural network is based on minimizing the loss function~--- the residual between the exact and approximate solutions for the training part of the dataset. These residuals can have different forms. Since the neural network ultimately solves an ill-posed inverse problem, the form of the loss function together with the dataset must ensure the uniqueness and stability of the solution. The form of the loss function also determines the quality (accuracy) of the neural network's operation. On the other hand, the loss function in our case must reflect the specificity of the \emph{jointness} of field inversion. Therefore, we will use so-called \textit{structural loss functions} to train the neural network. Let us say a few preliminary words about this.

In geophysics, for conventional (non-neural network) joint processing of heterogeneous input data, the so-called structural approach has proven itself well from the point of view of regularization properties (see, for example,~\cite{Haber97, gallardo03, bosch01}). An example of popular implementations of the structural approach is the use of the \emph{cross-gradients} functional, the ``\emph{Dice-Sorensen}'' coefficient~\cite{Dice45}, etc., as loss functions. These functionals define a measure of similarity of two sets and can be used to compare the spatial structures of gravitational and magnetic field sources with known substance densities and its magnetic characteristics. Important in the structural approach is the idea of using a composite residual that includes a number of functionals of different ``physical meaning'' with weights. This is what we will apply in our neural network algorithm.

Let us touch upon the problem of choosing the type and architecture of the neural network for solving the problem of joint inversion of potential fields. For solving multidimensional problems of separate inversion, convolutional neural networks (CNNs) have shown good results (see, for example,~\cite{lecun98,alexnet12}). The recent surge of interest in neural network data analysis is directly related to their application. In geophysics, data processing using CNNs is widely used for separate field inversion: seismic~\cite{lin18}, electromagnetic~\cite{noh19}, etc. For separate inversion of gravitational and magnetic fields, various CNN models have also been applied, and they have demonstrated their effectiveness (see \cite{he_cnn21,zhao25} and others). Of several architectures of such networks, the most promising, in our opinion, is the convolutional neural network with the abbreviation ``U-Net''~\cite{unet}, as well as various modifications of this architecture (Unet++~\cite{yufeng21}, Unet3++~\cite{huang20}). The application of such networks for solving three-dimensional geophysical exploration problems was described in~\cite{huang22}.

In our work, we propose an algorithm for joint inversion of gravitational and magnetic fields, based on the so-called \textit{two-level neural network approach}. The ``lower'' level includes two auxiliary neural networks of the U-Net type. The first computes the spatial distribution of gravitational field sources with a known constant substance density in a certain domain from measurements of this field on a remote surface. The other solves a similar problem for magnetic field sources with known constant magnetic characteristics in the same domain and from measurements on the same surface. Geometrically, the dataset for joint training of neural networks is built on the basis of various forms of synthetic ore bodies, generated using some stochastic algorithm from elementary objects (for example, cubes). For each such body, by solving forward problems, separate distributions of the gravitational and magnetic fields on the selected ``measurement surface'' are constructed and included in the dataset. The output results of the two neural networks, i.e. the source distributions, are combined for further training at the ``upper'' level. Specifically, they are used to compute the total loss function (weighted structural residual). The latter includes residual functionals, generally of different types, for gravitational and magnetic solutions, so as to ultimately ensure the proximity of the boundaries of the distributions of gravitational and magnetic field sources. Functionals reflecting various \textit{a priori} constraints can also be included in the loss function. This loss function is used to train the entire neural network system. Neural networks trained in this way can be applied to solving the inverse problem of finding an ore body with constant unknown substance density and magnetization.

The results of our work presented below are structured as follows. Section~\ref{sec:statement_of_ditect_problems} provides formulas for solving the forward problem for calculating gravitational and magnetic fields at given points from given distributions of gravitational and magnetic field sources. The general formulation of the problems for separately determining these sources from measurements of the corresponding fields in a remote region is considered. This involves not only finding their spatial distribution but also finding the variable masses of gravitational sources and the magnetization vectors of magnetic sources. An example is given showing that both separate and joint solution of these inverse problems has a non-unique solution. In this connection, the conditions of the inverse problem are refined. It is assumed that the masses of the sources are the same, as well as the magnitizations, and they (their estimate) are known in the first approximation. Section~\ref{sec:a_structural_two-level_neural_network} presents the details of the two-level neural network method for solving the inverse problem of joint field inversion using U-Net type neural networks. A description of the composition and method of constructing the dataset for training two independent neural networks that separately solve the problems of inversion of gravitational and magnetic fields (lower level of inversion) is given. The structural residual is defined, which, when solving the inverse problems jointly (at the upper level of solution), allows finding a single body and corresponding sources that create the observed fields. Section~\ref{sec:numerical_experiments} discusses the details of optimizing the operation of the used neural networks. Special attention is paid to various dataset optimization procedures. The results of model computational experiments using the proposed neural network method for joint inversion of potential fields are discussed. The results of joint field inversion for real geophysical data are also presented. Section~\ref{sec:numerical_experiments_for_sources_without_apriori_information} considers the formulations of problems of separate and joint field inversion in the case of unknown constant density and magnetization of the substance. Structural residuals in this case are described, and numerical experiments on solving corresponding model problems are conducted. Section~\ref{sec:discussion} is devoted to the discussion of the obtained results and conclusions.

\section{Formulation of forward problems of magnetometry and gravimetry} \label{sec:statement_of_ditect_problems}

The forward problem of computing the magnetic field from the source distribution is solved within the following mathematical model. We assume that some volume $D$ is filled with magnetic masses with magnetization intensity (magnetization vector) \[\textbf{\emph{M}}(\textbf{\emph{r}}) = \big(M_x(\textbf{\emph{r}}),  M_y(\textbf{\emph{r}}), M_z(\textbf{\emph{r}})\big)^T,\] where $\textbf{\emph{r}} = (x,y,z) \in D$. It is known (see, for example,~\cite{YBook}) that the magnetic induction $\textbf{\emph{B}}(\textbf{\emph{r}}_s) = \big(B_x(\textbf{\emph{r}}_s),  B_y(\textbf{\emph{r}}_s), B_z(\textbf{\emph{r}}_s)\big)^T$, created by these magnetic masses at a point $\textbf{\emph{r}}_s = (x_s,y_s,z_s) \notin D$, is represented as
\begin{equation}
	\label{ch1:field_problem_statement_M-vector}
	\textbf{\emph{B}}(\textbf{\emph{r}}_s) = \frac{\mu_0}{4\pi} \, \int\int\limits_D\int \bigg(\frac{3 \big(\textbf{\emph{M}}(\textbf{\emph{r}}), \textbf{\emph{r}} -\textbf{\emph{r}}_s\big)(\textbf{\emph{r}} -\textbf{\emph{r}}_s)}{|\textbf{\emph{r}} - \textbf{\emph{r}}_s|^5} - \frac{\textbf{\emph{M}}(\textbf{\emph{r}})}{|\textbf{\emph{r}} - \textbf{\emph{r}}_s|^3}\bigg)dv.
\end{equation}
Here $\mu_0$ is the known magnetic constant.

To compute the vector field $\textbf{\emph{B}}(x, y, z)$ at points $(x_{s_j},y_{s_j},z_{s_j})$, $j= \overline{1,S}$, which determine the positions of sensors measuring the magnetic induction $\textbf{\emph{B}} (x_s,y_s,z_s)$ created by the body, we discretize equation \eqref{ch1:field_problem_statement_M-vector} 
by the simplest procedure. The domain $D$ is divided into $N$ subdomains $D_i$, $i = \overline{1,N}$. The coordinates $(x_i, y_i, z_i)$ are the positions of some geometric center of the $i$-th subdomain, and its volume is $dv_i$. The values $\textbf{\emph{M}}(x_i,y_i,z_i)$ are specified at these centers. As a result, instead of equality~\eqref{ch1:field_problem_statement_M-vector}, we obtain a system of equalities
\begin{equation}\label{ch3:main_SLAE_M}
	\textbf{\emph{B}} (x_{s_j},y_{s_j},z_{s_j})=\frac{\mu_0}{4\pi} \sum\limits_{i=1}^N {\textbf K}^{(m)}(x_i, y_i, z_i, x_{s_j},y_{s_j},z_{s_j}) \, \textbf{\emph{M}}(x_i,y_i,z_i) \, dv_i ,\quad j= \overline{1,S},
\end{equation}
where
\begin{equation*}
	\begin{aligned}
		&{\textbf K}^{(m)}(x, y, z, x_s, y_s, z_s) = \frac{1}{r^5}\begin{bmatrix}3(x-x_s)^2-r^2 & 3(x-x_s)(y-y_s) & 3(x-x_s)(z-z_s) \\ 3(y-y_s)(x-x_s) & 3(y-y_s)^2-r^2 & 3(y-y_s)(z-z_s) \\ 3(z-z_s)(x-x_s) & 3(z-z_s)(y-y_s) & 3(z-z_s)^2-r^2\end{bmatrix}
	\end{aligned}
\end{equation*}
and
\begin{equation*}
	r = |\textbf{\emph{r}} - \textbf{\emph{r}}_s| = \sqrt{(x - x_s)^2 + (y - y_s)^2 + (z - z_s)^2}.
\end{equation*}

In what follows, we assume that the points $(x_{s_j},y_{s_j},z_{s_j})$ are located in a rectangle $\Pi$ on some fixed plane, remote from $D$ such that $ (x_{s_j},y_{s_j},z_{s_j}) \in \Pi =[X^{min}_s,X^{max}_s]\times [Y^{min}_s,Y^{max}_s]\times \{z_{s_j}=Z_s\}$, $D\cap\Pi=\varnothing$.

The formulas for the other forward problem~--- computing the gravitational field potential from the source distribution (see \cite{YBook})~--- are constructed in the same domain $D$ and the same rectangle $\Pi$ of sensor positions. If the domain $D$ is filled with masses with volumetric substance density $\bar{\rho}(\textbf{\emph{r}})$, $\textbf{\emph{r}} \in D$, then the potential $\bar{\varphi}(\textbf{\emph{r}}_s)$ of the gravitational field, induced by these masses at the sensor location point $\textbf{\emph{r}}_s = (x_s,y_s,z_s) \notin D$, is represented as
\begin{equation}
	\label{ch1:field_problem_statement_G}
	\bar{\varphi}(\textbf{\emph{r}}_s) = G \, \iiint\limits_D \dfrac{\bar{\rho}(\textbf{\emph{r}})}{|\textbf{\emph{r}} - \textbf{\emph{r}}_s|}dv.
\end{equation}
Here $G$ is the gravitational constant.

Using the same finite-dimensional approximation scheme for relation \eqref{ch1:field_problem_statement_G}, we obtain the following set of equalities:
\begin{equation}
	\label{ch3:main_SLAE_G}
	\bar{\varphi}(x_{s_j},y_{s_j},z_{s_j})=G \sum\limits_{i=1}^N K^{(g)}(x_i, y_i, z_i, x_{s_j},y_{s_j},z_{s_j}) \, \bar{\rho}(x_i,y_i,z_i) \, dv_i  , \quad j = \overline{1,S},
\end{equation}
where
\begin{equation*}
	\begin{aligned}
		&K^{(g)}(x, y, z, x_s, y_s, z_s) = \frac{1}{r} = \frac{1}{\sqrt{(x - x_s)^2 + (y - y_s)^2 + (z - z_s)^2}}.
	\end{aligned}
\end{equation*}

These formulas allow us to find potentials $\varphi_j =\bar{\varphi}(x_{s_j},y_{s_j},z_{s_j})$ and fields $\textbf{\emph{B}}_j = \textbf{\emph{B}}(x_{s_j},y_{s_j},z_{s_j})$ at sensor points $(x_{s_j},y_{s_j},z_{s_j})$, $j = \overline{1,S}$, for given values $\rho_i = \bar{\rho}(x_i, y_i, z_i)$, $\textbf{\emph{M}}_i = \textbf{\emph{M}}(x_i, y_i, z_i)$ and $dv_i$ in given set of points $(x_i, y_i, z_i)$, $i = \overline{1,N}$.

To formulate the inverse problems, we first consider relations~\eqref{ch1:field_problem_statement_M-vector} and~\eqref{ch1:field_problem_statement_G}. They can be viewed as equations for finding functions $\bar{\rho} (\textbf{\emph{r}}),\,\textbf{\emph{M}}(\textbf{\emph{r}})$ themselves from field measurements in the rectangle $\Pi$. However, simple examples show that the solution of this inverse problem is often not unique.

\textbf{Example.} Let the common support of sufficiently smooth functions $\rho (\textbf{\emph{r}})$ and $\textbf{\emph{M}}(\textbf{\emph{r}})$ be a ball of radius $a>0$ centered at the origin. We assume that these functions depend only on $r = |\textbf{\emph{r}}|$. Then for the gravitational potential $\varphi (r)$ the equation $\Delta \varphi (r) =  - 4\pi \rho (r)$ holds. Assuming the field $\textbf{\emph{B}}(r)$ is potential with potential $U(r)$, i.e. $\textbf{\emph{B}}(r)=\operatorname{grad}\, U(r)$, we write a similar equation: $\Delta U(r) =  - \mu_0 \operatorname{div} \textbf{\emph{M}}(r)$. These equations have solutions for $r > a$
\begin{equation*}
	\varphi (r)\, = \frac{{4\pi }}{r}\int\limits_0^a {\rho (t){t^2}dt}, \, \, U(r)\, = \frac{{{\mu _0}}}{r}\int\limits_0^a {\operatorname{div} \textbf{\emph{M}}(t)\,{t^2}dt}.
\end{equation*}
Obviously, the same potential $\varphi (r)$ can correspond to several different density distributions $\rho (r)$. Indeed, taking some distribution $\rho (r)$, we can add to it, without changing the potential, a function $\rho_0 (r)$ such that $\rho (r)+\rho_0 (r)>0$ and $\int\limits_0^a {{\rho _0}(t){t^2}dt}  = 0$. A similar ambiguity in determining the quantity $\operatorname{div} \textbf{\emph{M}}(r)$, and hence $\textbf{\emph{M}}(r)$, follows from the form of the solution for $U(r)$. As the supports of the functions $\rho (r)$ and $\operatorname{div} \textbf{\emph{M}}(r)$ is the same, the problem of joint field inversion will also be ambiguously solvable. A similar example of ambiguity can be given for the inverse problem of jointly solving the discrete equations~\eqref{ch3:main_SLAE_M} and~\eqref{ch3:main_SLAE_G}.

From the given example, it follows why these inverse problems in practice are solved on a priori given classes of functions $\rho (\textbf{\emph{r}})$, $\textbf{\emph{M}}(\textbf{\emph{r}})$. Often, these functions are assumed to be constant: $\rho (\textbf{\emph{r}}) = {\rho _0}$, $\textbf{\emph{M}}(\textbf{\emph{r}}) = {\textbf{\emph{M}}_0}$, since in this case there are uniqueness theorems for the solution of these inverse problems for some classes of bodies \cite{L1,L2,L3}. Below we will use this assumption in sections~\ref{sec:a_structural_two-level_neural_network} and~\ref{sec:numerical_experiments}, and therefore we will say a few words about the limits of its applicability. When processing data for real gravitational anomalies, the quantities $\rho_i=\rho_0=\mathrm{const}$ from formula \eqref{ch3:main_SLAE_G} are interpreted as the excess densities of the ore body elements over the density of the surrounding rocks. Thus, in the considered model, constants $\rho$ significantly different from zero can be found only for significant gravitational anomalies caused by a homogeneous ore body. Next, we indicate why we will assume that the vector $\textbf{\emph{M}}(x, y, z)=\textbf{\emph{M}}_0$, interpreted as the excess magnetization over the background one in the anomaly region, is constant. From geophysical studies, it is known that the magnetization of an ore body is induced by the global geomagnetic field of the Earth, which on the scale of typical geological exploration tasks ($\sim 10$ km) can be considered constant. As a consequence, all magnetization vectors $\textbf{\emph{M}}(x, y, z)$ at various points of the study area for such sizes of the studied body can be considered co-directional. This is what we will assume in Sections~\ref{sec:a_structural_two-level_neural_network} and~\ref{sec:numerical_experiments}, setting a constant direction for the unit vector of $\textbf{\emph{M}}$. This assumption allows specifying only the length ${m}_0$ of the vector $\textbf{\emph{M}}$. It will be associated, as data of the inverse problem, not with the entire vector function of the induced magnetic field $\textbf{\emph{B}}$, but only with its $z$--component, which we will denote as ${b}$: ${b}=B_z$. This model assumption makes it possible to work with gravitational and magnetic solutions of the same matrix size, which somewhat simplifies the form of the applied structural neural network approach.

\textbf{\textit{Remark.}} Separately, it should be noted that the assumptions made above imply that the studied ore bodies in the process of their formation did not experience destructive effects from tectonic processes, volcanism, and other similar processes~\cite{bulter1992}. That is, it is assumed that the studied ore bodies lack domains with different directions of residual magnetization, which, in the presence of significant deposits of ferromagnetic materials, can significantly distort the restored shape of the ore body. Thus, the used model has a limited scope of applicability. In particular, it may not be applicable in continental studies of rocks containing significant deposits of ferromagnetic materials. However, such studies are beyond the scope of this work.

So, below, in Sections~\ref{sec:a_structural_two-level_neural_network} and~\ref{sec:numerical_experiments}, we assume that $\rho (\textbf{\emph{r}}_{i})=\rho _{0}$, $\textbf{\emph{M}}(\textbf{\emph{r}}_{i})=m_{0}\textbf{\emph{n}}_{M}$, where $\rho _{0}$, $m_{0}$ are known constants and $\textbf{\emph{n}}_{M}$ is a known unit vector. In this connection, equalities~\eqref{ch3:main_SLAE_M} and~\eqref{ch3:main_SLAE_G} can be written in the form
\begin{equation}\label{Oper}
	\left\{
	\begin{aligned}
		&\rho _{0}A_{g}(D)=\varphi, \\
		&m_{0}A_{m}(D)=b,
	\end{aligned}
	\right.
\end{equation}
where the operator $A_{g}(D)$ computes the gravitational potential created by the combined body $D=\bigcup\limits_{i=1}^{N}D_{i}$ with $\rho =1$ at the sensor locations. Thus, the body is treated as a set of unit gravitational sources with positions $\textbf{\emph{r}}_{i}$. The operator $A_{m}(D)$ similarly computes at the sensor locations the values of the $z$--component of the magnetic induction vector function, created by the body $D$ with a unit vector $\textbf{\emph{M}}=\textbf{\emph{n}}_{M}$. Here the body is treated as a set of unit magnetic sources with their own positions $\textbf{\emph{r}}_{i}$. Then the problems of separate field inversion reduce to separately solving equations \eqref{Oper} with respect to $D$, i.e., with respect to the positions of the corresponding sources. In joint inversion, it is required to solve system \eqref{Oper} with respect to their common position.

Such problems of joint inversion of potential fields, based on models~\eqref{ch3:main_SLAE_M} and~\eqref{ch3:main_SLAE_G}, were previously solved using traditional methods in a number of works (see, for example,~\cite{Leonov-2025}). As an alternative to these methods, we apply in this work a method based on a \emph{two-level neural network algorithm}. We will apply this algorithm to two types of joint inversion problems. The first type is joint inversion for known constants $\rho _{0},m_{0}$. Having implemented this algorithm in the form of trained neural networks for such problems, we will then consider the question of refining the a priori constants $\rho _{0}$, $m_{0}$ when solving the joint inversion problem for unknown $\rho (\textbf{\emph{r}})$, $\textbf{\emph{M}}(\textbf{\emph{r}})$ using the same neural networks. It will be assumed that the limits of variation of these quantities are known.

\section{Two-level neural network method} \label{sec:a_structural_two-level_neural_network}

This section presents the details of a two-level neural network method for solving the inverse problem of joint field inversion using U-Net neural networks. At the lower level of the algorithm, we will train two independent U-Net neural networks that separately solve the gravimetric and magnetic problems from the system~\eqref{Oper}. We will describe the main characteristics of these networks.

\subsection{Dataset generation} The datasets for both neural networks are of the same type. Their generation is divided into two stages. In the first stage, a permissible geometric set of three-dimensional domains of the form $D=\bigcup\limits_{i=1}^{N}D_{i}$ is created, characterized by centers $\textbf{\emph{r}}_{i}$, volumes $dv_{i}$, and modeling the sought ore body. This process will be discussed separately below due to its non-triviality. Then, using the values $\textbf{\emph{r}}_{i}$, $dv_{i}$ and the data $\rho _{0},m_{0}$, formulas~\eqref{ch3:main_SLAE_G} and~\eqref{ch3:main_SLAE_M} allow us to compute the values of the gravitational or magnetic fields $\varphi_j = {\varphi}(x_{s_j},y_{s_j},z_{s_j})$ and $b_j = {b}(x_{s_j},y_{s_j},z_{s_j})$, $j = \overline{1,S}$, created by the body on the sensor grid. For brevity, these computed values will be denoted as $\varphi (D)$ and $b(D)$. Thus, to the set of domains of the geometric dataset $\{D\}$ we associate the set of forward problem solutions $\{\varphi (D), b(D)\}$. As a result, the gravimetric dataset $\mathcal{D}^{(g)}=\{ D,\varphi (D)\}$ and the magnetometric dataset $\mathcal{D}^{(m)}=\{D, b(D)\}$ are formed. Formally, they define the correspondence between the output and input data of the neural networks. They can be used for separate training of the two neural networks, gravitational and magnetic, capable of reconstructing the distributions of the corresponding field sources. In the second stage, these datasets can be combined into a new set $\mathcal{D}_{joint}=\{\mathcal{D}^{(g)}, \mathcal{D}^{(m)}\}$ for joint training of the networks.

Let us discuss the architecture of the used neural network and some features of selecting its parameters.

\subsection{Neural network architecture} As known from the universal approximation theorem~\cite{cybenko89}, even a single-layer network can approximate an arbitrary function with any predetermined accuracy. However, in the case of a single-level network, the number of parameters required to achieve the given accuracy may be so large that the task of training the neural network becomes extremely difficult from a computational point of view. As a consequence, the efficiency of using the neural network sharply decreases. In this regard, to reduce the number of network parameters, the number of neuron layers is increased in it, i.e., a transition is made to \textit{deep learning}. In many applications, such as image processing tasks, the input and output data of the network have a very complex form. The same applies to the problem of potential field inversion, where given fields and sought source have different spatial dimensions and complex spatial structure. In such cases, convolutional neural networks (CNNs) with deep learning are effectively used.

{As the basic architecture of the neural network, which must reconstruct the spatial structure of the ore body from the gravitational or magnetic field, we will use the ``U-Net'' type architecture~\cite{unet}. It is built according to an encoder-decoder scheme. First, the input information is processed by the encoder, which compresses the data into a \emph{latent representation} of lower dimensionality but with a larger number of layers. Operations of a special kind (convolutions, etc.) are performed on the data, after which the decoder brings the data from the compressed {latent} representation to the original dimensionality or close to the specified one, depending on the problem being solved.}

To implement the work of the encoder and decoder, it is necessary to convert all elements of the dataset into a special tensor form. The gravitational ${ \varphi}$ and magnetic ${ b}$ fields observed on the sensor grid are converted into tensors $\varphi_{tensor}$ and $b_{tensor}$ of dimensionality $1\times M^x_S\times M^y_S$, where $M^x_S \cdot M^y_S=S$) In turn, the distribution of the centers of gravitational $\textbf{\emph{r}}^{(g)}=\{r_{i}^{(g)}\}$ and/or magnetic field $\textbf{\emph{r}}^{(m)}=\{r_{i}^{(m)}\}$ sources, $i=\overline{1,N}$, in each body are converted into special tensors $r_{tensor}^{(g)}$ and $r_{tensor}^{(m)}$ of dimensionality $M^z_N\times M^x_N\times M^y_N$, so that $M^z_N M^x_N M^y_N=N $. The method of filling the indicated tensors with the corresponding grid values can be arbitrary, but it is more convenient to fill them in such a way that the ``slice'' of the tensor for each fixed value of the first index represents a set of grid values of the corresponding function at a fixed depth in the domain $D$.

\begin{figure}[t]
	\centering
	\includegraphics[width=0.95\textwidth]{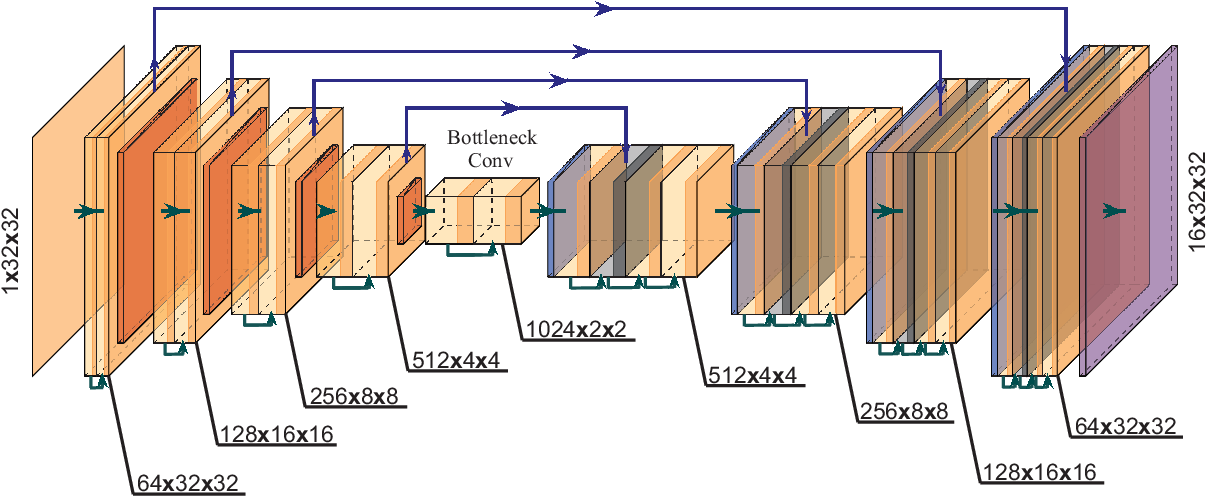}\\
	\caption{Architecture of the used ``U-Net'' type network by layers. The input is a tensor of dimensionality $1\times32\times32$. It corresponds to the grid values of size $32\times32$ of the gravitational or magnetic field measured in a rectangle, which is divided into cells by a two-dimensional Cartesian grid of dimensionality 32x32. The output of the network is a tensor of size $16\times32\times32$. It contains the grid values of the source density in a three-dimensional domain, which is divided into cells by a Cartesian grid of dimensionality $16\times32\times32$. }
	\label{unet}
\end{figure}

The datasets $\mathcal{D}^{(g)}=\{ D,\varphi (D)\}$, $\mathcal{D}^{(m)}=\{D, b(D)\}$ are converted to the form $\{r_{tensor}^{(g)},\varphi_{tensor}\}$,  $\{r_{tensor}^{(m)},b_{tensor}\}$ and are used separately for training the two neural networks. These networks implement the lower level of our field inversion algorithm. Their architecture is the same and is presented in Fig.~\ref{unet}. The procedure for transforming the dimensionality of data tensors in the encoder and decoder is also shown symbolically there. The principle of operation of encoders and decoders is described in more detail in the literature (see, for example, \cite{zafar22}, \cite[p. 233]{zhang21dive}).

\subsection{On training neural networks} At the lower level of the field inversion algorithm, two neural networks are trained for separate processing of gravitational and magnetic data. For this, the datasets $\mathcal{D}^{(g)}=\{ D,\varphi (D)\}$ and $\mathcal{D}^{(m)}=\{D, b(D)\}$ are used, each of which we symbolically represent as $\{ x,y\}$. Here $x$ is the exact solution of the inverse problem for its data $y$. Training consists of repeatedly minimizing some loss function $F(w) = F(\{ \mathord{\buildrel{\lower3pt\hbox{$\scriptscriptstyle\frown$}} \over x} (w)\} ,\{ x\} )$~--- the residual of the selected set $\{ x\}$ of exact solutions of the problem from the dataset and the set of approximate solutions $\{ \mathord{\buildrel{\lower3pt\hbox{$\scriptscriptstyle\frown$}} 	\over x} (w)\}$ computed by the neural network with parameters $w$ for the corresponding data $\{ y\}$ of the dataset. Minimization is performed over the network parameters $w$ sequentially for randomly selected sets $\{ x, y\}$ so that ultimately, when moving from set to set, the loss function decreases. Ideally, neural network parameters $\bar{w}$ should be obtained that can yield approximate solutions reproducing any element of the dataset. In practice, however, it is usually only possible to achieve an acceptable approximation accuracy,  that is a sufficiently small value of the loss function $F(w_{min})$ for the parameters $w_{min}$ of the trained network.

The so-called $\operatorname{Dice}$ functional is often used as the loss function. For a pair of vectors $\hat{x}=(\hat{x}_{1},...,\hat{x}_{L})$ and $x=(x_{1},...,x_{L})$ of the same dimensionality $L$, it is defined as
\begin{equation*}
	\operatorname{Dice}(\hat{x}, x) = \dfrac{2\sum
		\limits^L_{l=1} \hat{x}_l \, x_l}{\sum\limits^L_{l=1}(\hat{x}_l^2+ x_l^2)}.
\end{equation*}
Note that the axioms of a metric or norm are not satisfied for this functional.

We define a training datasets, i.e., some sets $\{ D_k^{(g)},\varphi (D_k^{(g)})\}$, $\{D_k^{(m)}, b(D_k^{(m)})\}$, $k=\overline{1,K}$. From a computational point of view, the sets $\{ D_k^{(g)}\}$, $\{D_k^{(m)}\}$ correspond to specifying groups of coordinates $\textbf{\emph{r}}_{k}^{(g)},~\textbf{\emph{r}}_{k}^{(m)}$ of the sources constituting the bodies $D_k^{(g)}$, $D_k^{(m)}$. Using the training samples, sets of source coordinates $\hat{\textbf{\emph{r}}}_{k}^{(g)}=\operatorname{NN}_{g}\big(\varphi (D_k^{(g)})\big)$, $\hat{\textbf{\emph{r}}}_{k}^{(m)}=\operatorname{NN}_{m}\big(b(D_k^{(m)})\big)$ can be calculated by the corresponding neural networks when applied to the data $\varphi (D_k^{(g)})$, $b(D_k^{(m)})$. These source coordinates correspond to the bodies $\hat{D}_k^{(g)}$, $\hat{D}_k^{(m)}$ found by the neural networks. Here and hereafter, the symbols $\operatorname{NN}_{g}$, $\operatorname{NN}_{m}$ denote the operators of the corresponding neural networks. Then the loss functions for our neural networks are defined as follows:
\begin{equation}\label{Dice}
	\begin{aligned}
		&\operatorname{Loss}_{grav}=\sum\limits_{k\in\text{\{batch\}}}
		\Big(1-\operatorname{Dice}\big(\hat{\textbf{\emph{r}}}^{(g)}_k, \textbf{\emph{r}}^{(g)}_k\big)\Big),
		&\operatorname{Loss}_{mag}=\sum\limits_{k\in\text{\{batch\}}}
		\Big(1-\operatorname{Dice}\big(\hat{\textbf{\emph{r}}}^{(m)}_k, \textbf{\emph{r}}^{(m)}_k\big)\Big).
	\end{aligned}
\end{equation}
Here the set of indices $\text{\{batch\}}$ is a subset of indices $k$ of fixed size from the set $\{\overline{1,K}\}$, which is determined randomly at each step (epoch) of neural network training.

For brevity, we will sometimes write formulas \eqref{Dice} in symbolic form
\begin{equation*}\label{Dice1}
	\begin{aligned}
		&\operatorname{Loss}_{grav}=\sum\limits_{k\in\text{\{batch\}}}
		\Big(1-\operatorname{Dice}\big(\hat{D}^{(g)}_k, D^{(g)}_k\big)\Big),
		&\operatorname{Loss}_{mag}=\sum\limits_{k\in\text{\{batch\}}}
		\Big(1-\operatorname{Dice}\big(\hat{D}^{(m)}_k, D^{(m)}_k\big)\Big),
	\end{aligned}
\end{equation*}
emphasizing thereby that the written residuals characterize the closeness of the geometries of the bodies from the training sample and the bodies found using the neural networks.

\subsection{Combining two neural networks and structural residuals}
During data inversion, the two independent neural networks yield two spatial distributions of sources: gravitational and magnetic, $\hat{\textbf{\emph{r}}}_{k}^{(g)}=\operatorname{NN}_{g}\big(\varphi (D_k^{(g)})\big)$ and $\hat{\textbf{\emph{r}}}_{k}^{(m)}=\operatorname{NN}_{m}\big(b(D_k^{(m)})\big)$. Geometrically, these distributions typically represent bodies that differ from each other. Our goal is to train both neural networks so that the shapes of these bodies are as close as possible to each other. For these purposes, we modify the residual.

Let us introduce a general loss function
\begin{equation}
	\label{loss_function_without_regularization}
	\operatorname{Loss} = \frac{1}{2}\operatorname{Loss}_{grav} +\frac{1}{2}\operatorname{Loss}_{mag}.
\end{equation}
It collectively characterizes how well both neural networks reproduce the model solutions of their inverse problems when working separately. However, this function does not reflect the degree of similarity of the shape (structure of the sets $\hat{\textbf{\emph{r}}}_{k}^{(g)}$ and $\hat{\textbf{\emph{r}}}_{k}^{(m)}$) of the reconstructed ore bodies for gravitational and magnetic data. This necessary similarity can be taken into account by adding an additional third term to~\eqref{loss_function_without_regularization}, resulting in a new function (``full structural residual'')
\begin{equation}
	\label{loss_function_with_regularization}
	\begin{aligned}
		\operatorname{Loss}_{\text{joint}} = &\sum\limits_{k\in\text{\{batch\}}}\Bigg(\frac{1}{2}\Big(1-\operatorname{Dice} \big(\hat{\textbf{\emph{r}}}_{k}^{(g)},{\textbf{\emph{r}}}_{k}^{(g)}\big)\Big)+ \frac{1}{2}\Big(1-\operatorname{Dice}\big(\hat{\textbf{\emph{r}}}_{k}^{(m)},{\textbf{\emph{r}}}_{k}^{(m)})+ \\&\hspace{+5.5cm}+\alpha\Big(1- \operatorname{Dice}\big(\hat{\textbf{\emph{r}}}_{k}^{(g)}, \hat{\textbf{\emph{r}}}_{k}^{(m)}\big)\Big)\Bigg).
	\end{aligned}
\end{equation}
The third term connects the ``gravitational'' and ``magnetic'' forms of the obtained ore bodies and ``penalizes'' their strong difference. This term ensures the structural connection between the gravitational and magnetic solutions. Symbolically, we write formula~\eqref{loss_function_with_regularization} as
\begin{equation*}\label{joint}
	\begin{aligned}
		\operatorname{Loss}_{\text{joint}} = &\sum\limits_{k\in\text{\{batch\}}}\Bigg(\frac{1}{2}\Big(1-\operatorname{Dice} \big(\widehat D_k^{(g)},D_k^{(g)}\big)\Big)+ \frac{1}{2}\Big(1-\operatorname{Dice}\big(\widehat D_k^{(m)},D_k^{(m)}\big)\Big)+ \\&\hspace{+5.5cm}+\alpha\Big(1- \operatorname{Dice}\big(\widehat D_k^{(g)},\widehat D_k^{(m)}\big)\Big)\Bigg).
	\end{aligned}
\end{equation*}

\begin{figure}[t]
	\centering
	\includegraphics[width=0.95\textwidth]{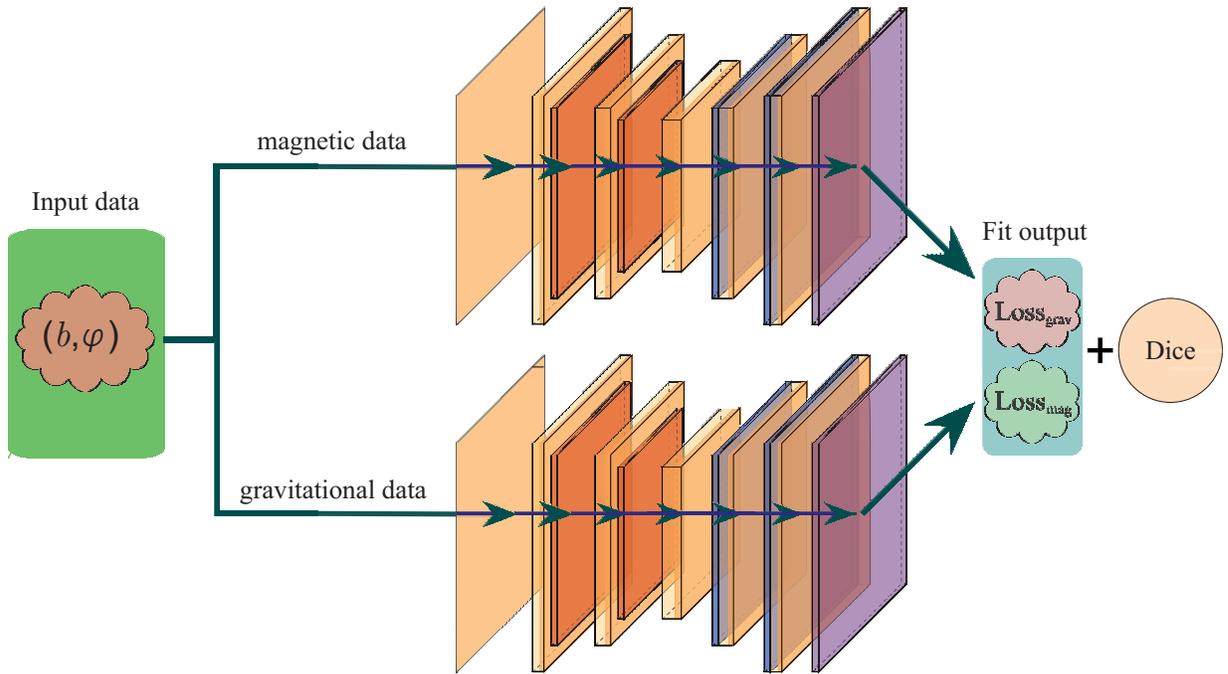}\\
	\caption{Training scheme of the proposed two-level neural network, including two ``U-Net'' type networks and a structural residual.}
	\label{unet_scheme}
\end{figure}

By training such an \emph{upper-level} network, we obtain a tool for joint field inversion. The architecture of our model with two neural networks and a common residual is presented in Fig.~\ref{unet_scheme}.
The parameter $\alpha>0$ is heuristic. Its choice is briefly discussed in Section~\ref{subsec:alpha_choice}.

\section{Numerical experiments and their analysis} \label{sec:numerical_experiments}

Let us detail the elements of the neural network algorithm used below for the numerical solution of model inverse problems of joint inversion of gravitational and magnetic fields. The neural networks used were trained for two joint inversion tasks. In the first task, the quantities $\rho_0$ and $m_0$ were assumed to be known and constant. It was required to determine the joint geometric position of the gravitational and magnetic field sources. In the second task, the quantities $\rho$ and $m$ were considered unknown constants lying within certain limits, which were to be found along with the joint geometric distribution of the sources. The second task was solved using the dataset generated for the first task. In Section~\ref{sec:numerical_experiments_for_sources_without_apriori_information} we present a discussion of this problem separately for the gravitational field. Now let's discuss the construction of a dataset to solve both of the mentioned problems.

\subsection{Dataset construction details}

In training neural networks to determine the three-dimensional shape of an ore body during joint field inversion, one of the key points is the selection/construction of the dataset. The ideal solution should be to include data on real ore bodies and solutions of the corresponding inverse problems. However, due to the lack of this information, it is necessary to generate and use synthetic datasets. The geometric shapes of the bodies included in it should reflect theoretically possible variants for which the corresponding potential fields can be computed analytically. This can be implemented (see~\cite{YBook}) by combining prisms of various types and sizes in various combinations of mutual arrangement, bodies in the form of steps, which are also combined with prisms, etc. In our work, a stochastic process similar to that presented in~\cite{huang21} is used for this purpose.

Another method, convenient to use with a limited set of training data and small resources, is \textit{transfer learning}. Its idea is as follows. A neural network is trained on some dataset that may not exactly correspond to the true data of the considered inverse problem. Then, to achieve sufficient solution accuracy, the network is further trained on a small set of data close to those that will be used in reality. This approach is widely used in geophysics for seismic problems~\cite{siahkoohi19}. However, for inverse problems of gravimetry and magnetometry, it is not sufficiently studied, and our work analyzes this approach for various types of data.

To describe the used datasets, let us define the domains in which the field sources and sensors are located. All geometric dimensions are given in meters. The other used quantities are dimensionless. The domain $V=[X^{min},X^{max}]\times[Y^{min},Y^{max}]\times[Z^{min},Z^{max}]$, in which the sought ore body lies, is specified by the quantities: $X^{min}=0$, $X^{max}=1\,600$, $Y^{min}=0$, $Y^{max}=1\,600$, $Z^{min}=0$, $Z^{max}=8\,000$. This domain was divided into $32\times32\times16=16\,384$ cells. The corresponding output tensor of the neural network has dimensions $M^z_N=16$, $M^x_N=32$, $M^y_N=32$. For the model problems being solved, each cell is a cube with a side of 50. The body $D$ itself is modeled by various unions of cubes (see below) with constant values $\rho_0$, $m_0$. For cubes lying outside the body the values $\rho=0$, $m=0$ are assigned.

The rectangle $\Pi= [X^{min}_s,X^{max}_s] \times[Y^{min}_s,Y^{max}_s]\times\{ z_{s_j} = Z_s\}$, on which the model sensors are located, was chosen with $X_s^{min}=0$, $X_s^{max}=1\,600$, $Y_s^{min}=0$, $Y_s^{max}=1\,600$, $Z_s=0.1$. It was divided into $32\times32 =1\,024$ cells, and this corresponds to the input tensor of the neural network with dimensions $M^z_N=1$, $M^x_N=32$, $M^y_N=32$.

Several types of datasets were used.
Ore bodies from the {TOY} dataset were generated geometrically according to the following algorithm:
\begin{enumerate}
	\item In the domain $V$, coordinates of one or two ``centers'' were randomly set. In the vicinity of each center, 4 cube-cells of size 2x2x2 were generated with the same unit value of source density in all cells.
	\item Each cube was iteratively (40 iterations) shifted in a randomly chosen direction by 2 cells.
\end{enumerate}

The body obtained as a result of these actions was included in the geometric dataset. A total of $K=11\,000$ such randomly formed 3D objects were built. Using the formulas from Section~\ref{sec:statement_of_ditect_problems}, the gravitational and magnetic fields corresponding to these bodies were computed and included in the datasets $\mathcal{D}^{(g)}=\{ D_k^{(g)},\varphi (D_k^{(g)})\}$, $\mathcal{D}^{(m)}=\{D_k^{(m)}, b(D_k^{(m)})\}$, $k=\overline{1,K}$. The {TOY} dataset was used to solve the first task.

In addition, more primitive geometric datasets were used. In particular, the {SYN} dataset consisted of a sufficiently representative set of combinations of prisms and steps. Figures close to elements from {SYN}, but obtained by stochastic combination of a small number of these elements, constituted the {STOCH} dataset. These datasets were used below for fine-tuning neural networks during the implementation of the \textit{transfer learning} process (see subsection~\ref{subsec:transfer_learning}).

\subsection{Neural network training details}

The well-known optimizer $AdamW$ (see, for example,~\cite{loshchilov18}) was used in optimizing the parameters $w$ of the neural network, with a convergence rate coefficient $learning\, rate=3.0 \cdot 10^{-4}$. The dataset $\mathcal{D}$ was split into training and test sets: $\mathcal{D} =\mathcal{D}_{\text{train}}\cup \mathcal{D}_{\text{test}}$, with $10\,000$ and $1\,000$ objects in each, respectively. Thus, the loss functions $\operatorname{Loss}_{\text{train}}$ and $\operatorname{Loss}_{\text{test}}$ like~\eqref{joint} with $\mathcal{D}_{\text{train}}$ and $ \mathcal{D}_{\text{test}}$ where used.

The size of the set $\{batch\}$ in these formulas was set to 64, the number of training epochs was set as $\mathrm{epoch}=300$. In this case, the best residuals were obtained in numerical experiments. At each epoch of neural network training, the corresponding final (minimum) values $\operatorname{Loss}_{\text{train}} = \operatorname{Loss}_{\text{train}}(\text{epoch}) $ and $\operatorname{Loss}_{\text{test}} = \operatorname{Loss}_{\text{test}}(\text{epoch})$ are computed.

\begin{figure}[t]
	\centering
	\includegraphics[width=0.95\textwidth]{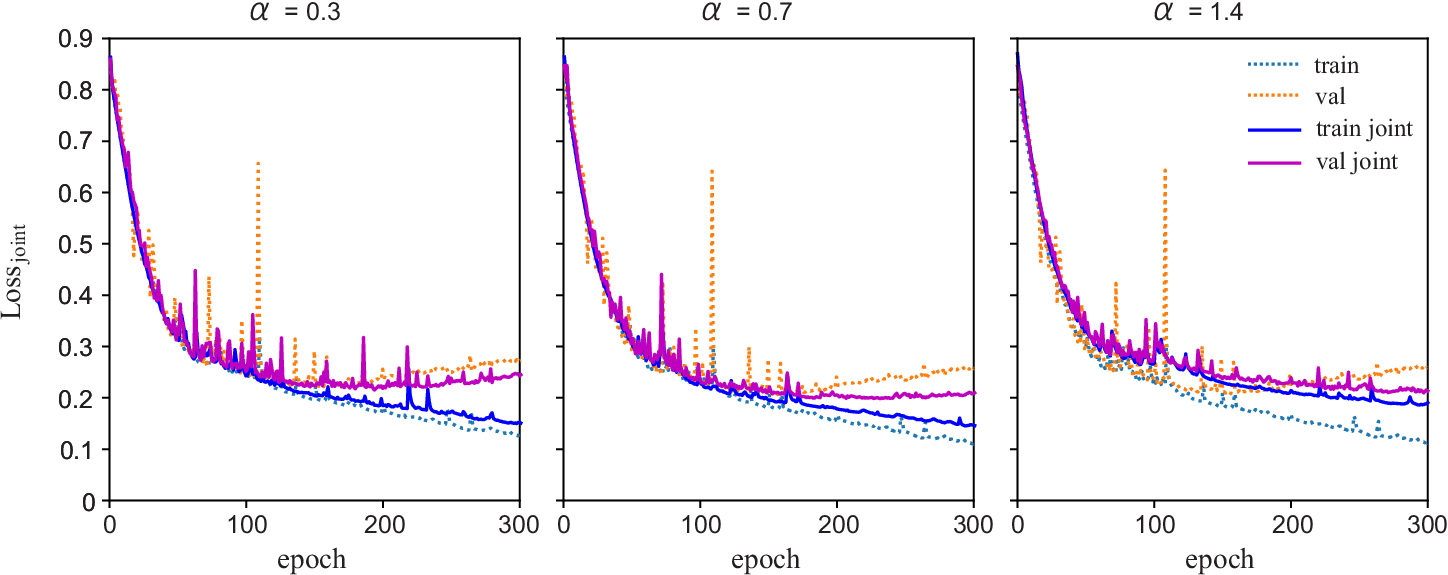}\\
	\caption{Comparison of the behavior of the curves $\operatorname{Loss}_{\text{train}}(\text{epoch})$ and $\operatorname{Loss}_{\text{test}}(\text{epoch})$ for different values of the parameter $\alpha$ (solid lines) and for $\alpha=0$ (dots).}
	\label{loss_dice}
\end{figure}

When determining the optimal moment to stop training the neural network, it is important to capture the situation when the algorithm begins to \emph{overfit} the training data. As a rule, this moment is determined by the beginning of the ``divergence'' of the curves $\operatorname{Loss}_{\text{train}}(\text{epoch})$ and $\operatorname{Loss}_{\text{test}}(\text{epoch})$ (see Fig.~\ref{loss_dice}). In our algorithm, training was stopped when the inequality $\big|\operatorname{Loss}_{\text{train}}(\text{epoch}) - \operatorname{Loss}_{\text{test}}(\text{epoch})\big| \geqslant \varepsilon$ was first satisfied. The number $\text{epoch}_{epsilon}$ found in this way determines the end of the iterations: $\operatorname{iter}_{\text{stop}} = \text{epoch}_{epsilon}$. Taking into account the form~\eqref{loss_function_with_regularization} of the loss function, it can be understood that the characteristic values of $\operatorname{Loss}_{\text{train}}$ and $\operatorname{Loss}_{\text{test}}$ do not exceed 1. Therefore, $\varepsilon=0.02$ was chosen, which corresponds to a relative difference in $\operatorname{Loss}$ residuals of a few percent.

The final residual value, determining the quality of the approximate solution, was taken as $\operatorname{Loss}_{\text{result}}=\operatorname{Loss}_{\text{test}}(\operatorname{iter}_{\text{stop}})$.

\subsection{Model example of reconstructing sources determining the three-dimensional structure of an ore body} \label{subsec:alpha_choice}

\begin{figure}[t]
	\centering
	\includegraphics[width=0.95\textwidth]{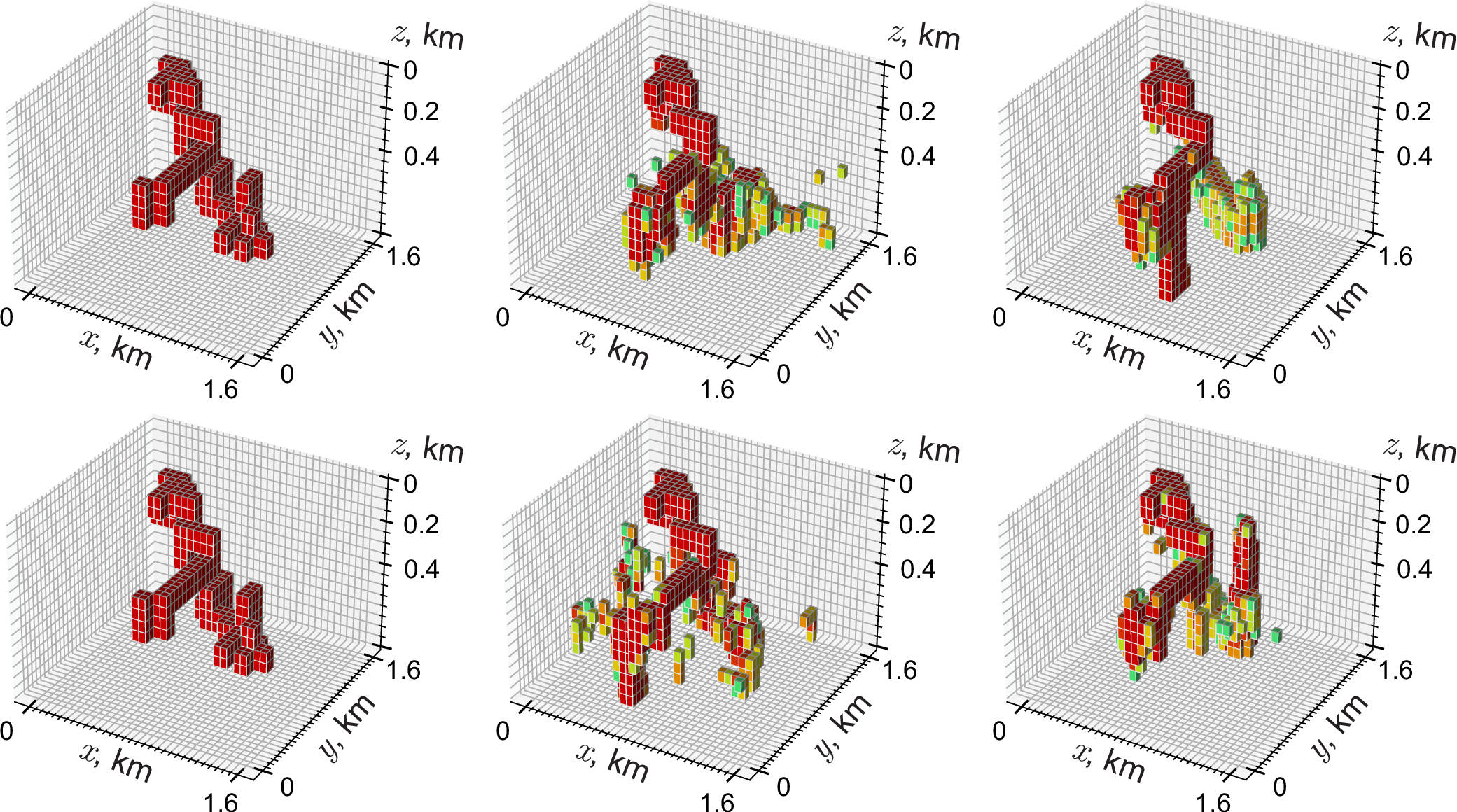}\\
	\caption{An example of reconstructing gravimetric and magnetic sources. The upper figures correspond to the distribution of gravitational field sources, the lower ones represents the distribution of magnetic field sources. The left column shows the model exact solutions, the middle column shows the result of separate source reconstruction, and the right column shows the result of joint reconstruction.}
	\label{example}
\end{figure}

Let us consider the reconstruction of the gravitational and magnetic field sources determining the object shown in the left column of Fig.~\ref{example}. For this object, we compute the input data for the neural network, i.e., the gravitational and magnetic fields in the rectangle $\Pi$, and first process them separately using the corresponding network, and then jointly, using the trained neural network with the structural residual~\eqref{loss_function_with_regularization} at $\alpha = 1$. The results are presented in columns 2 and 3 of Fig.~\ref{example}. It can be seen that in the case of separate solution of the inverse problems, the magnetic field sources are reconstructed with higher accuracy due to the better conditioning of the magnetometric problem (column 2). This is reflected in the accuracy of the corresponding neural network. The overall spatial structure of the model ore body is better reconstructed during joint data processing using the two-level network (the third column of the figures).

The calculation results significantly depend on the form of the structural norm, and especially on the contribution of its term with the parameter $\alpha$. Therefore, let us return to Fig.~\ref{loss_dice}, which illustrates the influence of this parameter on the residuals $\operatorname{Loss}_{\text{train}}$ and $\operatorname{Loss}_{\text{test}}$ of the full neural network. With an increase in the parameter $\alpha$, a noticeably smaller divergence of the curves $\operatorname{Loss}_{\text{train}}(\text{epoch})$ and $\operatorname{Loss}_{\text{test}}(\text{epoch})$ is observed, indicating an increase in the stability of the results. At the same time, the iteration number $\operatorname{iter}_{\text{stop}}$, marking the divergence of the curves $\operatorname{Loss}_{\text{train}}(\text{epoch})$ and $\operatorname{Loss}_{\text{test}}(\text{epoch})$, increases. Thus, better accuracy will be achieved at sufficiently ``large'' values of the parameter $\alpha$. However, this also increases the computational complexity of training the neural network.

\begin{figure}[t]
	\centering
	\includegraphics[width=0.75\textwidth]{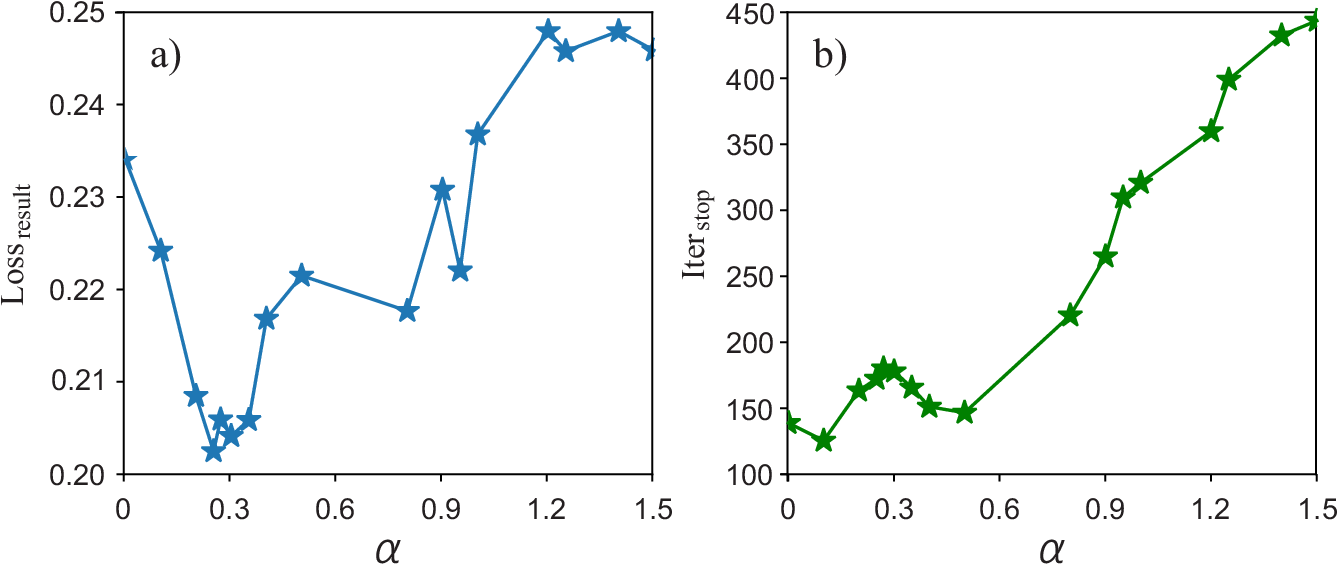}
	\caption{Left: dependence of accuracy $\operatorname{Loss}_{\text{result}}$ on $\alpha$; Right: dependence of the iteration number $\operatorname{iter}_{\text{stop}}$ on $\alpha$.}
	\label{err_alpha_dice}
\end{figure}

We optimized the parameter $\alpha$ heuristically. For this, the dependence of the neural network solution accuracy $\operatorname{Loss}_{\text{result}}$ on the parameter $\alpha$ was computed on a certain grid of values $\alpha \in [0,1.5]$, and then the best accuracy was selected. In Fig.~\ref{err_alpha_dice} on the left, the values of the best accuracy $\operatorname{Loss}_{\text{result}}=\operatorname{Loss}_{\text{test}}(\operatorname{iter}_{\text{stop}})$ for different values of $\alpha$ are presented. As can be seen, the reduction in the residual at the minimum point of this curve compared to the residual for $\alpha=0$ is about 8\%. Thus, the introduction of the structural connection (the third term) in the residual~\eqref{joint} provides a significant improvement in the solutions of the joint neural network. In Fig.~\ref{err_alpha_dice} on the right, the dependence of the final training iteration number $\operatorname{iter}_{\text{stop}}$ on $\alpha$ is shown. This number on average increases with $\alpha$, confirming the increasing complexity of training.

\subsection{Experiments with transfer learning} \label{subsec:transfer_learning}

We conducted initial experiments about the impact of transfer learning on the validation accuracy of the training dataset only for the lower-level neural network solving the inverse gravity problem. The experimental design was as follows. The neural network was trained on various mixtures of the STOCH and TOY datasets described above, each with the same number of elements, $K$. The composition of the mixture was formed in the percentage ratio $\lambda\in [0,1]$, i.e. it consisted of $\lambda\cdot K$ elements of the STOCH set and $(1-\lambda)\cdot K$ elements of the TOY dataset. Thus, the resulting dataset $\mathcal{D}(\lambda)$ contains geometric solutions $D_{k}^{(g,STOCH)}$, $D_{k}^{(g,TOY)}$ from the STOCH and TOY datasets mixed in the proportion $\lambda$ and the corresponding mixed gravitational field data. We denote by ${D}_{k}^{(g)}(\lambda )$, $k=1,\ldots,K$, the geometric elements of the mixed dataset $\mathcal{D}(\lambda)$, and $\varphi_{k}(\lambda )$ the corresponding gravity data.

The dataset $\mathcal{D}(\lambda)$ was split into two sets of indices $k$: the training set $\left\{learning\right\} $ and the validation set $\left\{valid\right\} $. For each $\lambda\in [0,1]$ defined on some grid, the neural network was trained on the $\left\{learning\right\}$ part of the $\mathcal{D}(\lambda)$ dataset. Then, the trained neural network was validated. The elements $D_{k}^{(g,STOCH)}$, $D_{k}^{(g,TOY)}$ of the STOCH and TOY sets for $k\in \left\{ valid\right\}$ were selected as exact geometric solutions for validation. They were compared with the corresponding solutions $\hat{D}_{k}^{(g)}(\lambda )$ at index $k\in \left\{ valid\right\}$ obtained by the neural network for the corresponding mixed data $\varphi_{k}(\lambda )$. The residual
\begin{equation*}
	Loss \equiv
	Loss_{valid}^{(test)}(\lambda )=\sum\limits_{k\in \{valid\}}\Big( 1-Dice\big( \hat{D}_{k}^{(g)}(\lambda ),D_{k}^{(g,test)}\big) \Big),
\end{equation*}
which determines the accuracy of validation, was chosen as the comparison measure. Here the symbol $test$ is either STOCH or TOY. The choice of different validation groups $\{D_{k}^{(g,STOCH)}\}$, $\{D_{k}^{(g,TOY)}\}$, $k\in \left\{ valid\right\}$, allows us to compare the influence of the mixture coefficient $\lambda$ on the validation accuracy for different datasets, i.e. to evaluate the influence of $\lambda$ on the accuracy of geometric solutions of a certain class. For comparison, a third validation group of elements from the SYN dataset, which is often found in other works, was also used.

\begin{figure}[t]
	\centering
	\includegraphics[width=0.75\textwidth]{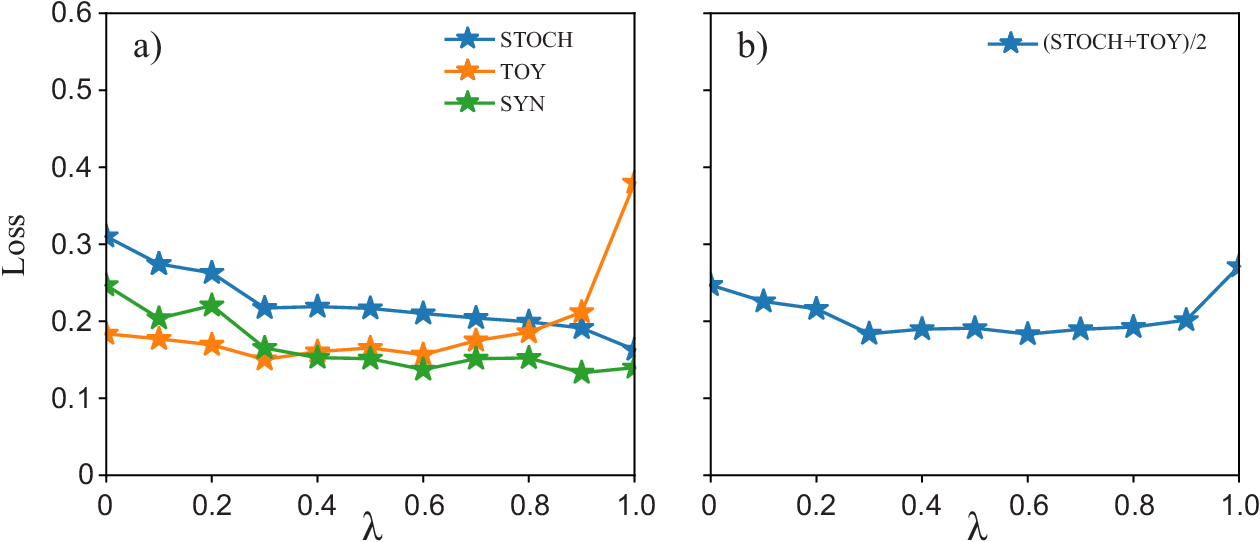}
	\caption{Left: influence of the percentage ratio $\lambda$ on the validation accuracy (the $\mathrm{Loss}$ function). Right: averaged validation value for mixed validarion dataset (50 \% STOCH + 50 \% TOY).}
	\label{Fig2}
\end{figure}

The results of validation on the STOCH, TOY and SYN classes are shown in Fig.~\ref{Fig2} on the left. The graphs show that, in general, validation on the TOY class provides better accuracy than on the STOCH class. This is expected, since the STOCH class consists of more primitive objects than the TOY class. The effect is also visible of a sharp decrease in validation accuracy in the TOY validation at $\lambda = 0.9$ when the training dataset contains a relatively large number of elements from STOCH (see the line TOY). At the same time, an increase in the volume of TOY in the training set, when $\lambda = 0.3$, improves the accuracy of validation. This illustrates the application of the transfer learning approach, when adding thirty percent of the STOCH elements to the TOY data in the training set improves validation accuracy. Note that the total volume of the training set does not change, i.e. the same volume of resources is used for training in all cases.

Validation on the SYN class turns out to be better than the others when $\lambda>0.4$. This is explained by the fact that SYN objects have a simpler form, which is better reconstructed by the neural network. In a certain sense, the SYN set can be considered a subclass of the STOCH set. Therefore, after training on a sample with a relatively high content of STOCH elements ($\lambda>0.5$), the neural network is better at reconstructing SYN elements during validation than TOY elements.

Figure~\ref{Fig2} on the right shows the results of validation on another set of elements, consisting of 50 percent of elements from STOCH and TOY. For values of $\lambda\in [0,0.8]$, this validation accuracy is slightly worse than the accuracy of validation using TOY. The best accuracy is achieved at $\lambda=0.6$.

\subsection{Influence of including noised data in the dataset on solution accuracy}

Another way to improve a dataset to increase validation accuracy is to include noisy data. Numerical experiments of this type were also conducted for the lower-level neural network solving the inverse gravity problem. Let us introduce a new type of dataset, STOCH NS. It contains geometric objects from STOCH and their corresponding fields. Furthermore, the same geometric objects are associated with the same fields perturbed by normally distributed noise with a zero mean and a root-mean-square error of 0.02.

\begin{figure}[t]
	\centering
	\includegraphics[width=0.75\textwidth]{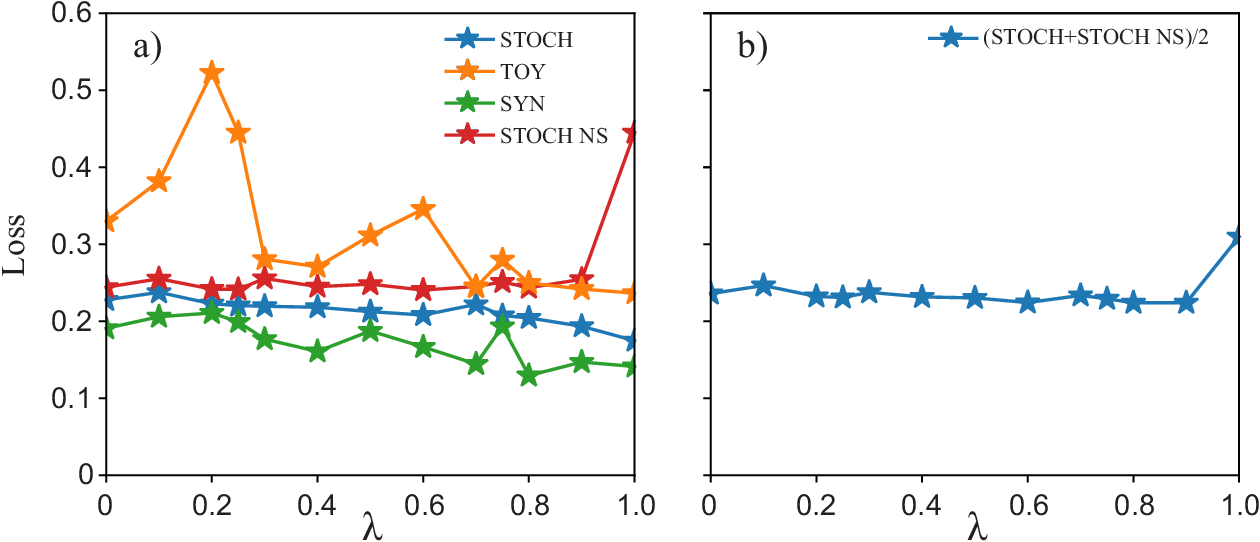}
	\caption{Left: influence of the percentage ratio $\lambda$ in the STOCH NS dataset on validation accuracy. Right: averaged validation value for the two datasets STOCH, STOCH NS.}
	\label{Fig3}
\end{figure}

For training, mixed dataset of volume $11\,000$ was formed, consisting of elements from STOCH and STOCH NS in a percentage ratio $\lambda$. For $\lambda=0$ it consists only of STOCH NS elements, and for $\lambda=1$ it containes only elements from STOCH. As in the previous section, after training, for each $\lambda$, the neural network was validated on sets from STOCH, TOY, and SYN, as well as from STOCH NS. The validation results are presented in Fig.~\ref{Fig3}.

In the left figure, it can be seen that the validation accuracies on the noised STOCH NS and non-noised STOCH turn out to be close even at a small percentage of STOCH NS ($\lambda\in (0,0.9)$) in the training dataset. Validation on the STOCH dataset is best with small additions of STOCH NS elements to the training dataset. Conversely, a significant increase in the number of STOCH NS elements in the training dataset generally worsens validation accuracy on both the STOCH and SYN datasets. Thus, the number of perturbed elements in the training dataset should not be large. The right side of Fig.~\ref{Fig3} shows the results of validation on a dataset consisting of 50 percent elements of the STOCH and STOCH NS sets.

\subsection{Experimental data processing}\label{Brazil}

\begin{figure}[h]
	\centering
	\includegraphics[width=0.75\textwidth]{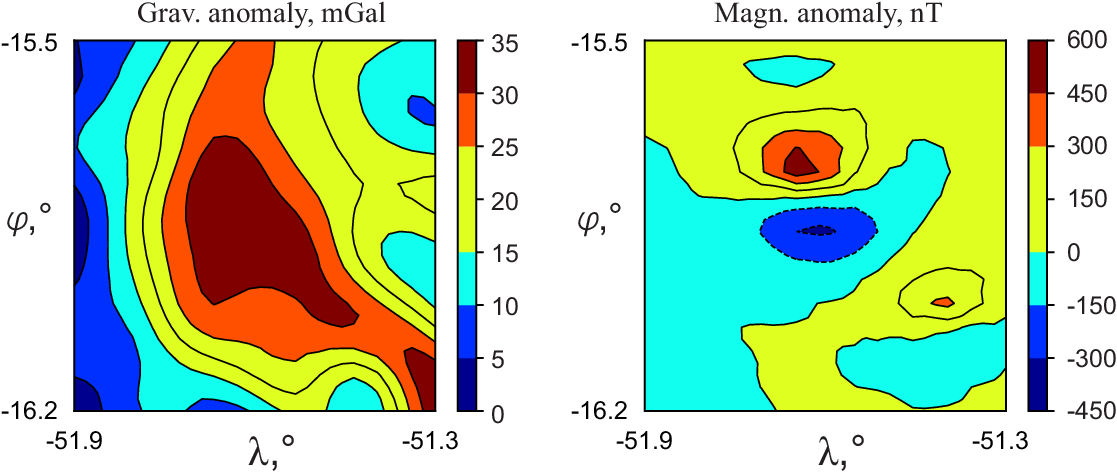}
	\caption{Data of the inverse problems for the Jussara region, Goias State, Brazil: measurements of the $z$--components of the fields.}
	\label{Real_data}
\end{figure}

As an example of the application of extended datasets with noised data, consider the problem of practical joint inversion of gravitational and magnetic fields when processing data for the region of Jussara, Goias State, Brazil with latitude and longitude $\varphi\in [-16.2^o,-15.5^o]S, \lambda\in [-51.9^o, -51.25^o]E$. The data for these problems~--- measurements of the $z$-components of the fields~--- are given in Fig.~\ref{Real_data}. The data and grids were taken from the WGM2012 GLOBAL MODEL [International Gravimetric Bureau. https://bgi.obs-mip.fr/data-products/grids-and-models/wgm2012-global-model] (gravity data) and WDMAM [Wdmam. http://wdmam.org/ 

\noindent download.php] (magnetic data) databases.

The conversion of gravitational data, i.e., the $z$--components $g_z(x_s,y_s,z_s)$ of the gravitational field in the sensor domain $\Pi$, into potential values $\varphi(x_s,y_s,z_s)$ required for neural network processing, was carried out using the known formula
\begin{eqnarray*}
	\varphi (x_{s},y_{s},z_{s})=-\frac{1}{2\pi }\int \int_{\Pi }\frac{g_{z}(x^{\prime },y^{\prime },z_{s})dx^{\prime }dy^{\prime }}{\sqrt{(x_{s}-x^{\prime })^{2}+(y_{s}-y^{\prime })^{2}}}\approx -\frac{1}{2\pi }\sum\limits_{k=1}^{K}C_{k}(x_{s},y_{s})g_{z}(x_{k},y_{k},z_{s}).
\end{eqnarray*}
Here $C_{k}(x_{s},y_{s})=\int \int_{\Pi _{k}}\frac{dx^{\prime }dy^{\prime }}{\sqrt{(x_{s}-x^{\prime })^{2}+(y_{s}-y^{\prime })^{2}}}$, and $\Pi _{k}$ are rectangles with centers $(x_{k},y_{k})$, constituting a sufficiently detailed partition of the domain $\Pi =\bigcup\limits_{k=1}^{K}\Pi _{k}$.

To jointly solve inverse problems, we used two differently trained neural networks described in previous subsections, NN and NNns with noisy data. The neural networks $\operatorname{NN}$ were trained on a sample of size $22\,000$ objects from the {STOCH} dataset with input field data dimensions of $32\times 32$ and solution dimension of $16\times32\times32$. The function~\eqref{Dice} were used as the loss function. When applying the neural network $\operatorname{NNns}$ to process real data, the sample size was $66\,000$ objects. The proportion of noised objects was $0.1$ of the total. The noise level was set at $2\%$. The residual for this dataset is used according to~\eqref{Dice}.

\begin{figure}[h]
	\centering
	\includegraphics[width=0.65\textwidth]{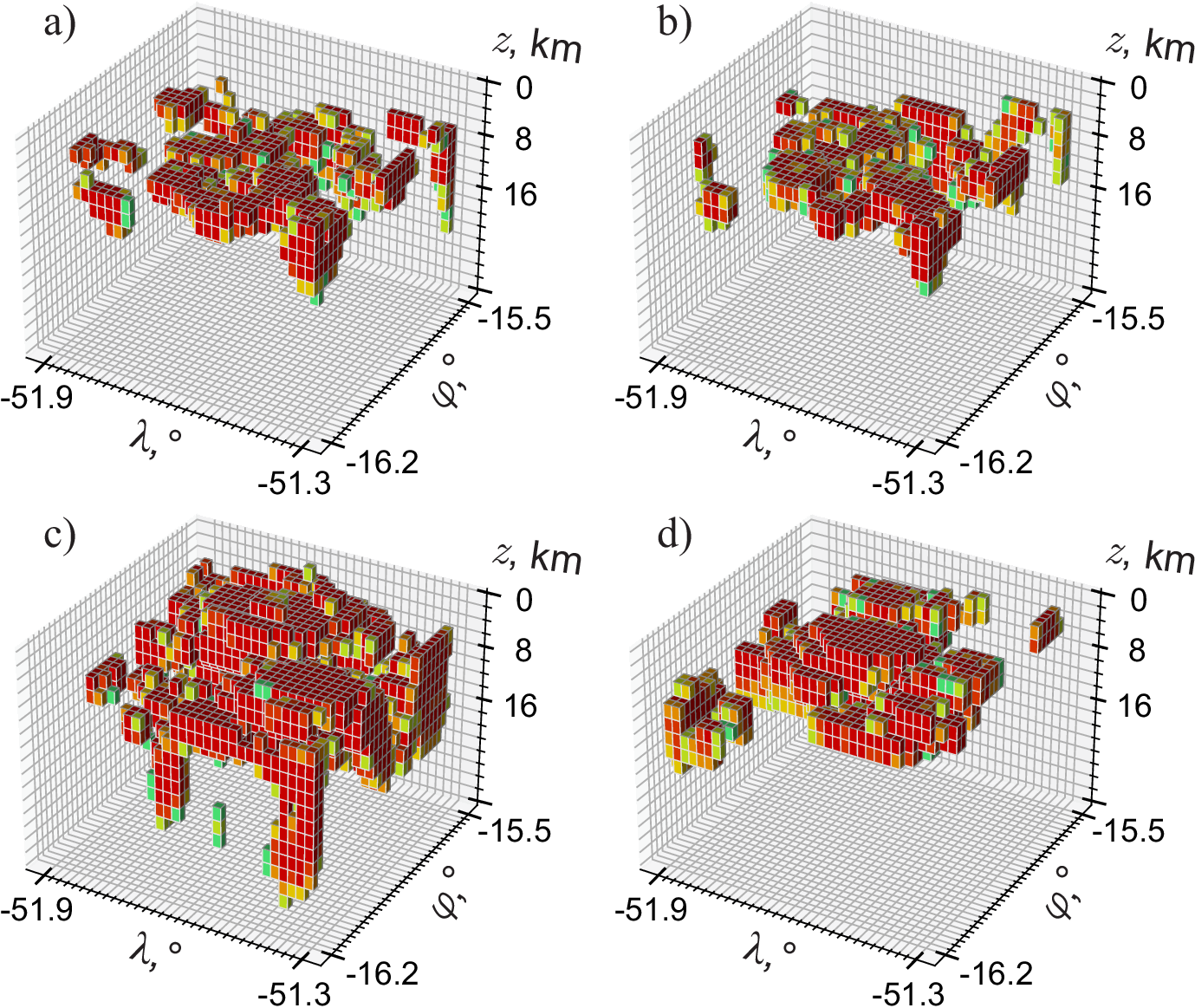}
	\caption{Processing of data for Jussara, Goias State, Brazil. The top row shows gravimetric results of joint field inversion using the neural network $\operatorname{NN}$ (first column) and the neural network $\operatorname{NN}ns$ (second column). The bottom row shows the same for magnetometric results of joint inversion.}
	\label{Fig6}
\end{figure}

Fig.~\ref{Fig6} presents the results of the joint determination of the shape of the ore body (at $\rho=1$) from this data using these neural networks.

A significant improvement in quality is noticeable for the neural network $\operatorname{NNns}$ for gravitational data. For magnetic data, the results of the neural networks $\operatorname{NNns}$ and $\operatorname{NN}$ differ insignificantly.

\section{Numerical experiments for problems with unknown $\rho$, $m$}  \label{sec:numerical_experiments_for_sources_without_apriori_information}

In this section, we consider the inverse problem of jointly determining gravitational and magnetic sources with unknown positions and unknown constants $\rho$ and $m$. This corresponds to searching for the geometry of the ore body $D$, as well as its gravitational and magnetic properties $\rho$ and $m$ from the data $\varphi$ and $b$ measured in the domain $\Pi$. Mathematically, this is the joint solution with respect to $D$, $\rho$, and $m$ of the operator equations
\begin{equation*}\label{unknown_dencity}
	\rho A_g (D)=\varphi, \quad m A_m (D)=b.
\end{equation*}
Here, interval constraints of the form $\rho _{0}-\Delta \rho \leq \rho \leq \rho _{0}+\Delta \rho$, $m_{0}-\Delta m\leq m\leq m_{0}+\Delta m$ are imposed on the unknowns $\rho$ and $m$, with given $\Delta \rho ,\Delta m$. To solve this problem, the neural networks from Section~\ref{sec:numerical_experiments} will be used.

\subsection{On separate estimating source density}

As an example, let us first consider the inverse gravimetry problem~\eqref{ch3:main_SLAE_G} for the case of a \emph{constant}, but unknown, value $\rho$ of the excess density of the ore body over the surrounding rocks. We assume that the ore body is represented by one of the elements $D$ of some set of three-dimensional domains $\mathcal{D} = \left\{ {D}\right\}$. Then problem~\eqref{ch3:main_SLAE_G} for each fixed $\rho$ can be written in a form similar to~\eqref{Oper}: find $D_\rho^*\in \mathcal{D}$ such that
\begin{equation}\label{Azu}
	{D_\rho^*} = \operatorname{argmin} \big\{ d (\rho A_g[D],\varphi):\,\,D \in \mathcal{D} \big\}.
\end{equation}
Here, as in formulas~\eqref{Oper}, $A_g[{D}]$ is the operator computing the distribution of the gravitational potential created by the body ${D}$ of unit density in the observation domain, $\varphi$ is the inverse problem data, and $d(\cdot,\cdot)$ is the used measure of proximity (residual) between the data and their computed analogs. Problem~\eqref{Azu} is interpreted as follows: find the body (bodies) $D_\rho^*\in \mathcal{D}$ with constant density $\rho$ that creates in the observation domain $\Pi$ a field $\rho A_g[D]$ that best approximates (in the sense of the residual $d(\cdot,\cdot)$) the data $\varphi$. Due to the ambiguity of simultaneously determining both the value $\rho$ and the domain ${D}$ in problem~\eqref{ch3:main_SLAE_G} (see the example from Section~\ref{sec:statement_of_ditect_problems}), problem~\eqref{Azu}, generally speaking, will have different solutions $D^*_\rho\in \mathcal{D}$ for different $\rho$. The same effect occurs when using neural networks to solve problem~\eqref{ch3:main_SLAE_G} in the form~\eqref{Azu}. Let us illustrate this with numerical experiments.

\begin{figure}[t]
	\centering
	\includegraphics[width=0.65\textwidth]{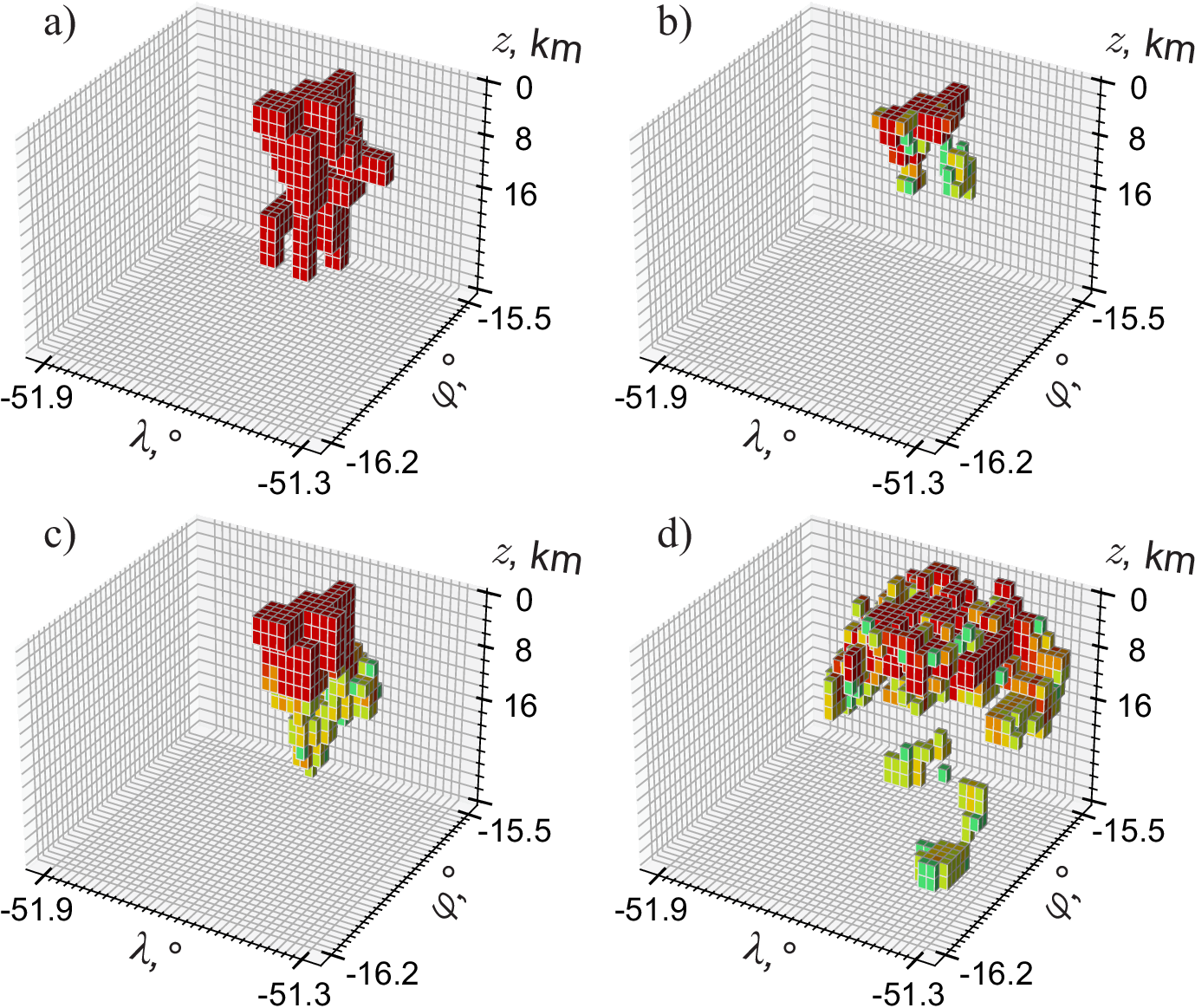}
	\caption{Change in the geometry ${D_\nu^*}$ of sources as the value of $\nu$ changes: a) model solution, b) $\nu = 3$, c) $\nu = 1$, d) $\nu = \frac{1}{3}$.}
	\label{Fig4}
\end{figure}

Suppose a neural network $D=\operatorname{NN}_g(\varphi)$ is trained to solve problem \eqref{Azu} with known  $\rho=\rho_0$ on the dataset $\{\mathcal{D},\varphi(\mathcal{D})\}$. Let us apply it to solve this problem on the same dataset, but with other constant densities $\rho=\nu\cdot\rho_0>0$ for various $\nu>0$. Here we assume that the values $\nu$ vary within known limits such that the constraints $\rho _{0}-\Delta \rho \leq \rho \leq \rho _{0}+\Delta \rho$ are satisfied. This means that we are looking for a solution $D^*_\nu\in \mathcal{D}$ of the extremum problem
\begin{equation}\label{Azu1}
	{D_\nu^*} = \operatorname{argmin} \Big\{ d\big(\nu\rho_0 A[\operatorname{NN}(\varphi)],\varphi\big) \Big\}.
\end{equation}
for various $\nu$. The results of these experiments are shown in Fig.~\ref{Fig4}.
It is clear from it that the geometric distribution of sources $D$ found by the neural network nonlinearly depends on $\nu$, and as $\nu$ increases, the domain ${D_\nu^*}$ decreases in volume.

Now suppose the data $\tilde{\varphi}$ are not taken from the training set, i.e., $\tilde{\varphi}\notin \varphi(\mathcal{D})$.
Then one can pose the question of selecting the value of the average source density of the form $\rho=\nu\rho_0$ that ensures the best solution of problem \eqref{Azu1} with data $\tilde{\varphi}$ processed by the neural network $\operatorname{NN}_g$, which was trained for $\rho=\rho_0$. The corresponding problem consists in minimizing the function
\begin{equation*}
	\bar{\Phi}_g (\nu ) = d\big( {\nu\rho_0A_g[\operatorname{NN}_g(\tilde{\varphi})], \tilde{\varphi} }\big)
\end{equation*}
within the limits $\nu  \in [{\nu _1},{\nu _2}]$, produced by the constraints on $\rho$. In practice, the limits $\nu_1$, $\nu_2$ for the variation of $\nu$ are determined by geological data about the underlying rocks.

Similar considerations hold for the inverse magnetometry problem included in~\eqref{unknown_dencity}. Therefore, introducing $m=\mu m_{0}$, one should minimize the function
\begin{equation*}
	\bar{\Phi}_m (\mu ) = d\big({\mu m_0A_m[\operatorname{NN}_m(\tilde{b})], \tilde{b}}\big)
\end{equation*}
to find the unknown value $\mu$ using the neural network $\operatorname{NN}_m$ within the known limits $\mu \in \lbrack \mu_{1},\mu_{2}]$. Here, in general, $\tilde{b}\notin b(D)$.

We will consider two types of residuals $d(\cdot ,\cdot )$ in parallel: $d_1(b,\breve{{b}})=\| b-\breve{{b}}\|_{2}^{2}$ and $d_2(\varphi,\breve{\varphi})=1-\mathrm{Dice}(\varphi,\breve{\varphi})$, and use the functions $\bar{\Phi}^{(1,2)}_g (\nu ),\bar{\Phi}^{(1,2)}_m (\mu )$, generated by these residuals respectively.

\begin{figure}[t]
	\centering
	\includegraphics[width=0.65\textwidth]{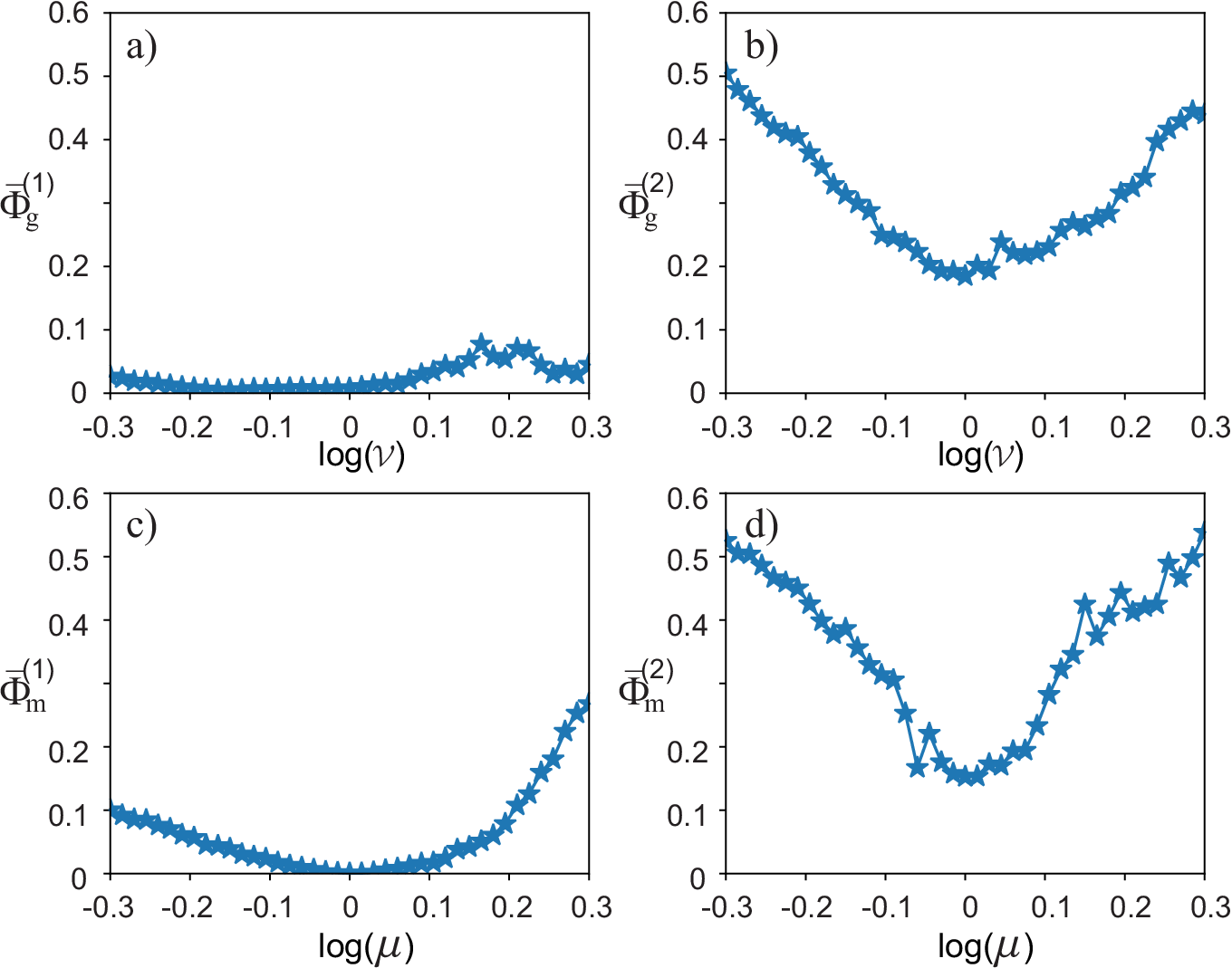}
	\caption{Averaged over the validation set residuals $\bar{\Phi}^{(1,2)}_g (\nu )$ and $\bar{\Phi}^{(1,2)}_m (\mu )$ for the gravitational and magnetic problems.}
	\label{fig:psi_phi_resid}
\end{figure}

The next group of experiments was as follows. The same model problem of field inversion as in subsection~\ref{subsec:alpha_choice}  with constants $\rho_0=1$ and $m_0=1$ was solved, and a common body $D_0$ containing gravitational and magnetic sources was found using neural networks of the type $\operatorname{NNns}$. Then the values $\rho =\nu \rho _{0}$, $m=\mu m_{0}$ were refined by minimizing the functions
\begin{equation*}
	\bar{\Phi}_g (\nu ) = d\big( {\nu\rho_0A_g[D_0], \tilde{\varphi} }\big), \quad
	\bar{\Phi}_m (\mu ) = d\big( {\mu m_0A_m[D_0], \tilde{b} }\big)
\end{equation*}
for $(\nu ,\mu )\in \lbrack \nu _{1},\nu _{2}]\times \lbrack \mu _{1},\mu _{2}]$.
In Fig.~\ref{fig:psi_phi_resid}, the left column shows the residuals $\bar{\Phi}^{(1)}_g (\nu )$ and $\bar{\Phi}^{(1)}_m (\mu )$ averaged over the validation set, depending on the values of $\nu$, $\mu$. The right column shows similar values of the functionals $\bar{\Phi}^{(2)}_g (\nu )$ and $\bar{\Phi}^{(2)}_m (\mu )$.
From the figures, it is clear that with this approach, the density $\rho=1$, i.e., $\nu=1$, is better recovered using the function $\bar{\Phi}^{(2)}_g (\nu )$. For the magnetic problem, both functions $\bar{\Phi}^{(1,2)}_m (\mu )$ provide good recovery of the value $m=1$, i.e., $\mu=1$.

Let us emphasize that in the presented calculations, the problem of non-uniqueness in simultaneously finding the density $\rho$, magnetization $m$ of the ore body, and its geometry is, to a certain extent, solved by fixing the possible body forms, i.e., through the geometric datasets of the neural networks. \textit{A priori} information on the values of $rho_0$ and $mu_0$ was also significantly used.

\begin{figure}
	\centering
	\includegraphics[width=0.65\textwidth]{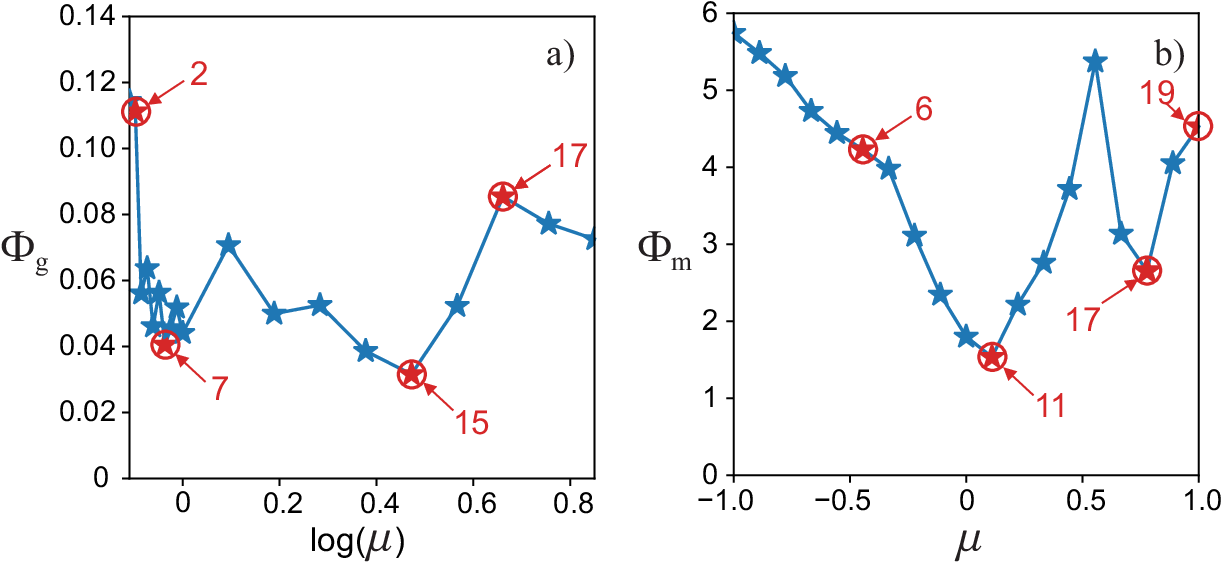}
	\caption{Averaged over the validation set residuals $\bar{\Phi}^{(1)}_g (\nu )$ and $\bar{\Phi}^{(1)}_m (\mu )$ for the gravitational (a) and magnetic (b) problems. }
	\label{Brasil_separate}
\end{figure}

\begin{figure}[p]
	\centering
	\includegraphics[width=0.65\textwidth]{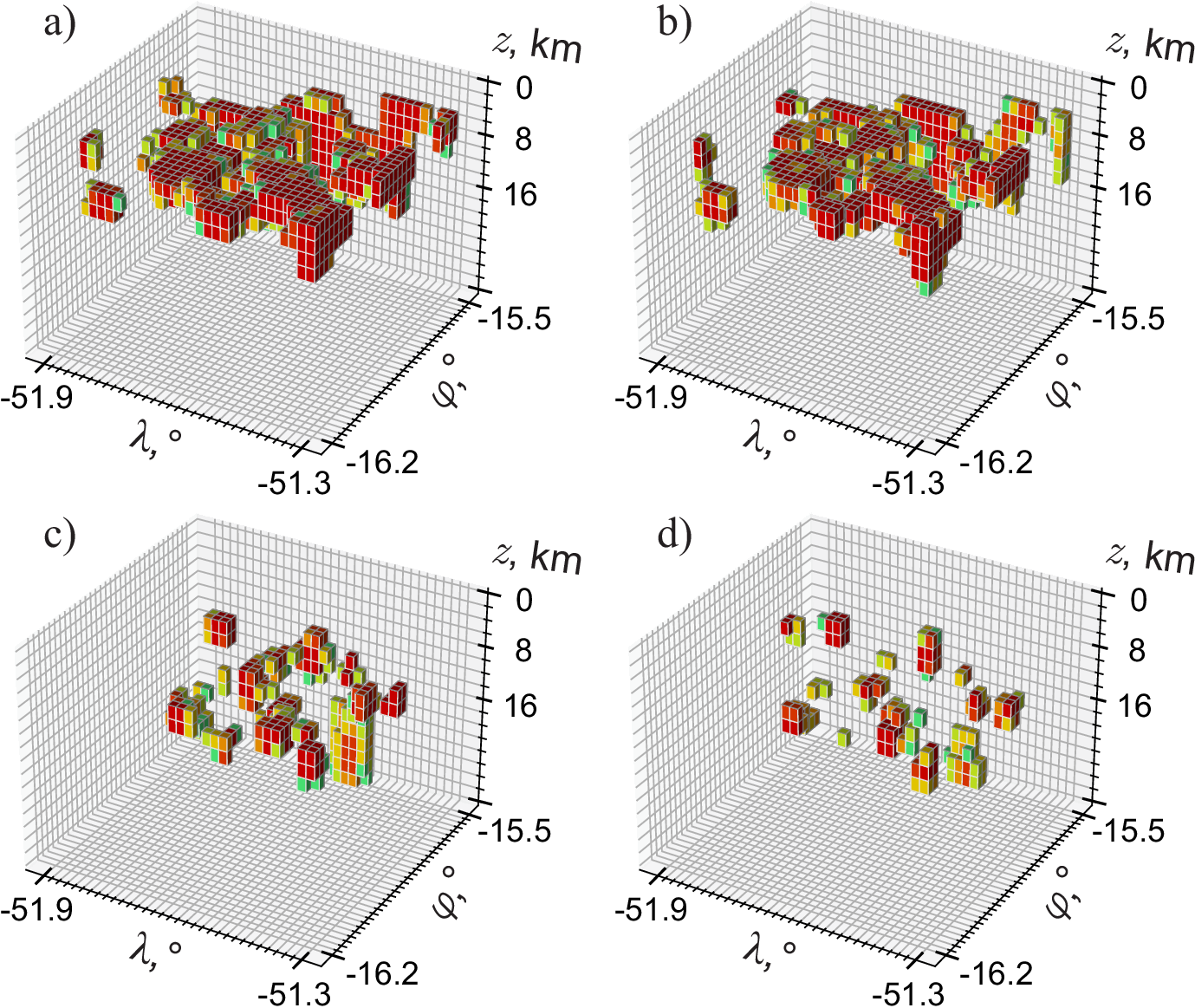}
	\caption{Gravimetric forms of ore bodies for selected values of $\nu$. These solutions correspond to the following highlighted points in Fig.~\ref{Brasil_separate}: a) 2, b) 7, c) 15, d) 17. }
	\label{G_domains}
\end{figure}

\begin{figure}[p]
	\centering
	\includegraphics[width=0.65\textwidth]{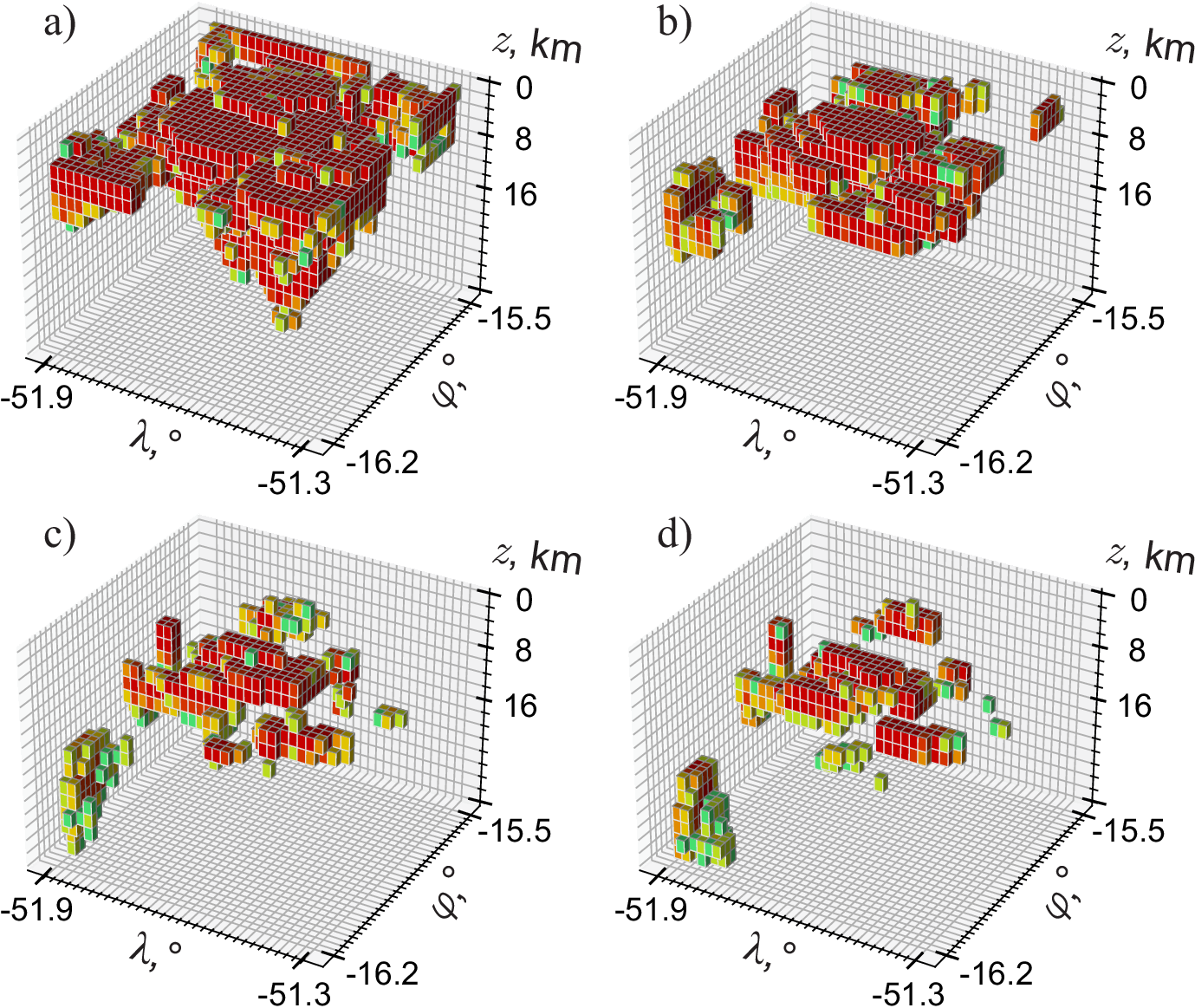}
	\caption{Magnetometric forms of ore bodies for selected values of $\mu$. These solutions correspond to the following highlighted points in Fig.~\ref{Brasil_separate}: a) 6, b) 11, c) 17, d) 19.}
	\label{M_domains}
\end{figure}

The above methodologies for simultaneously finding the field sources and the values $\rho$ and $m$ were applied for the \emph{separate processing} of experimental data from the Jussara region, Goias State, Brazil (see subsection~\ref{Brazil}). The averaged residuals $\bar{\Phi}^{(1)}_g (\nu )$ and $\bar{\Phi}^{(1)}_m (\mu )$ computed from these data for $rho_0=1$ and $mu_0=1$ are shown in Fig.~\ref{Brasil_separate}.
The figures show several local minima of these functions, which is explained by the non-uniqueness of the solution for each inverse problem. Let us consider the shapes of the ore bodies obtained for some values of $\nu$ and $\mu$. In the case of the magnetic field (see Fig.~\ref{Brasil_separate}-b), a minimum stands out for the point with \emph{number 11}. The shape of the ore body for it is shown in Fig.~\ref{M_domains}-c. In the gravitational case (Fig.~\ref{Brasil_separate}a)), the global minimum of the function is at point 15, which corresponds to the body in Fig.~\ref{G_domains}-c. There is also a local minimum at point 7 with a set of oscillations, which corresponds to the ore body in Fig.~\ref{G_domains}-b. It is most similar in shape to the bodies from the magnetic problem. The difference in the shapes of the obtained bodies is associated with the separate solution of the field inversion problems. Fig.~\ref{G_domains},~\ref{M_domains} also show other shapes of bodies corresponding to local minima from the graphs in Fig.~\ref{Brasil_separate}.

\subsection{Joint inversion of gravimetric and magnetic data with unknown $\rho$ and $m$.}

Now let us consider problem~\eqref{unknown_dencity} of \emph{joint inversion} of fields with unknowns the average densities $\rho$ and $m$ of the gravitational and magnetic sources. Here we will use a priori information in the form of estimates $\rho_0$ and $m_0$ for the densities of the gravitational and magnetic sources. In problems~\eqref{unknown_dencity}, we introduce again quantities $\nu$, $\mu$ such that $\rho =\nu \rho _{0}$, $m=\mu m_{0}$. Formally, the problem of joint field inversion~\eqref{unknown_dencity} when searching for the body shape $D$ and the values $\rho$, $m$ using neural networks of the type $\operatorname{NN}_g(\varphi )$, $\operatorname{NN}_m(b )$ and $\operatorname{NNns}_g(\varphi )$, $\operatorname{NNns}_m(b )$ reduces, as in the previous subsection, to some joint minimization of the functions $\Phi_{g}(\nu)=\bar{\Phi}^{(1)}_g (\nu )$, $\Phi_{m}(\mu)=\bar{\Phi}^{(1)}_m (\mu )$ or $\Phi_{g}(\nu)=\bar{\Phi}^{(2)}_g (\nu )$, $\Phi_{m}(\mu)=\bar{\Phi}^{(2)}_m (\mu )$. This can be achieved by constructing new composite residual functions of various types. They should include the functions $\bar{\Phi}^{(1,2)}_g (\nu )$ and $\bar{\Phi}^{(1,2)}_m (\mu )$ themselves, for example in the form
\begin{equation*}
	\Phi_{full}(\nu,\mu)=\Phi_{g}(\nu)+\Phi_{m}(\mu)
\end{equation*}
and minimize it with respect to the parameters $\mu$, $\nu$. To ensure the proximity of the geometry of the bodies obtained when processing gravitational and magnetic data, functions of the form
\begin{equation*}
	S(\nu,\mu) = d(\operatorname{NN}_g(\tilde \varphi), \operatorname{NN}_m(\tilde b)), \quad
	S(\nu,\mu) = d(\operatorname{NNns}_g(\tilde \varphi), \operatorname{NNns}_m(\tilde b)) \label{eq:struct_reg}
\end{equation*}
can be included in the composite residual.
Then joint inversion can be carried out by minimizing a composite function of the type
\begin{equation}\label{nev_joint}
	\Phi_{joint}(\nu,\mu) = \beta_1\Phi_{g}(\nu) + \beta_2\Phi_{m}(\mu) +\alpha S(\nu,\mu)
\end{equation}
for $(\nu ,\mu )\in \lbrack \nu _{1},\nu _{2}]\times \lbrack \mu _{1},\mu _{2}]$. Here $\beta_1,\beta_2,\alpha>0$ are weight coefficients characterizing the contribution of each of the residuals included in $\Phi_{joint}(\nu,\mu)$. The choice of these weights is usually done heuristically.

\begin{figure}[t]
	\centering
	\includegraphics[width=0.85\linewidth]{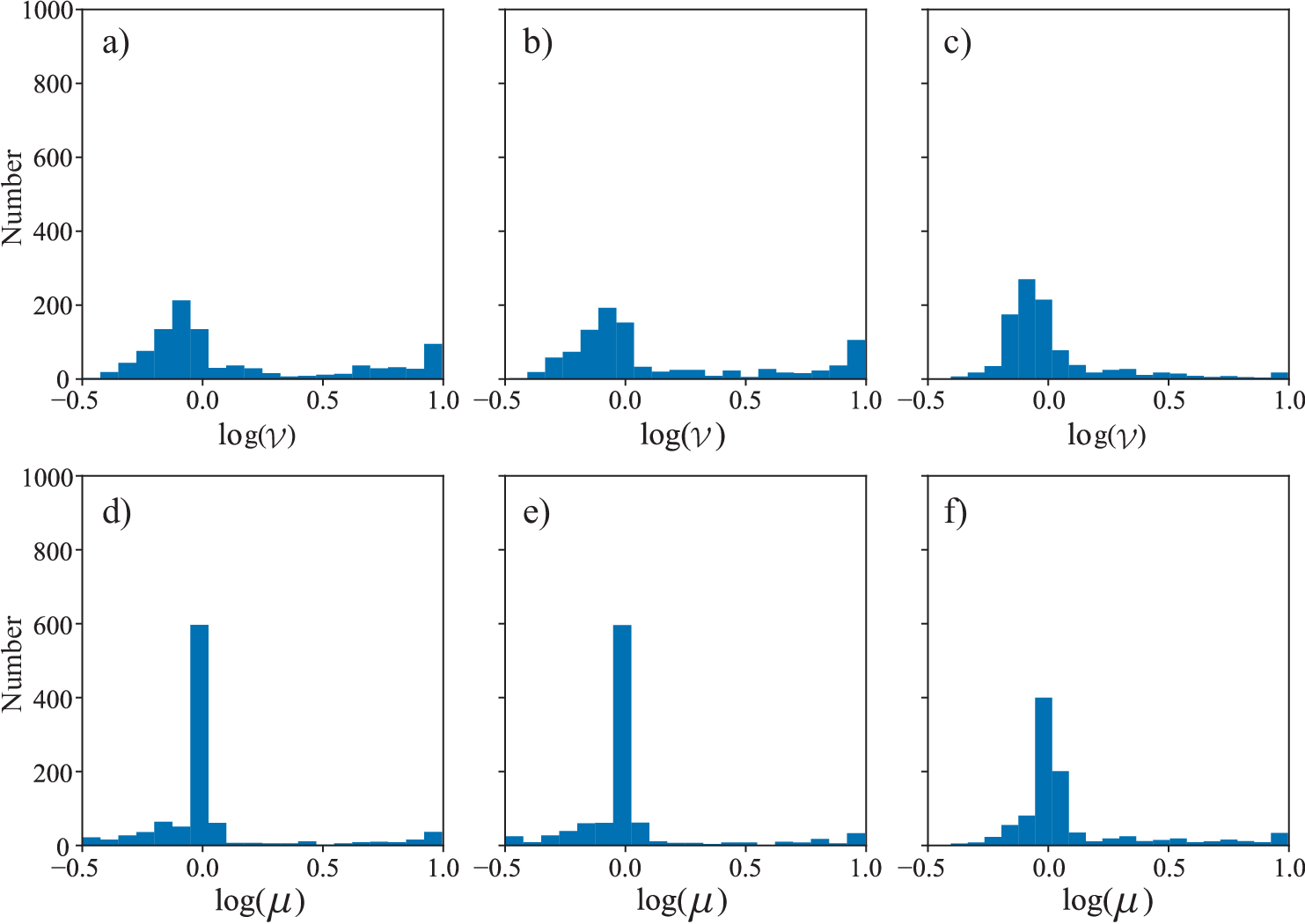}
	\caption{Histograms of the recovery of values $\nu$, $\mu$ when using various functions~\eqref{nev_joint}.  }
	\label{fig:histo_joint}
\end{figure}

Again, we consider as an example the model problem from subsection~\ref{Brazil}. By repeatedly solving the minimization problem of function \eqref{nev_joint} for various data $\{\varphi ,b\}$ on grids of values $\nu ,\mu $ under the constraints $-0.5\leq \nu$, $\mu \leq 1$, we obtain the distribution of $\nu$ and $\mu$ values in the form of histograms in Fig.~\ref{fig:histo_joint}. The histograms show the empirical frequency of occurrence of certain values of the parameters $\nu$, $\mu$. The first row of figures relates to the processing of gravitational data, and the second to magnetic data. The columns of the figure contain results for different sets of weights $\beta_1$, $\beta_2$, $\alpha$. Thus, the first column corresponds to the case of using in the residual \eqref{nev_joint} only the functions $\Phi_{g}(\nu)=\bar{\Phi}^{(1)}_g (\nu )$, $\Phi_{m}(\mu)=\bar{\Phi}^{(1)}_m (\mu )$ ($\beta_1=\beta_2=1,\alpha=0$), the second represents only the functions $\Phi_{g}(\nu)=\bar{\Phi}^{(2)}_g (\nu )$, $\Phi_{m}(\mu)=\bar{\Phi}^{(2)}_m (\mu )$ ($\beta_1=\beta_2=1,\alpha=0$), and the third is associated with the use of three functions $\Phi_{g}(\nu)=\bar{\Phi}^{(2)}_g (\nu )$, $\Phi_{m}(\mu)=\bar{\Phi}^{(2)}_m (\mu )$, $S(\nu,\mu)$ with $\beta _{1}=\beta _{2}=1$, $\alpha =0.2$.
It can be seen from Fig.~\ref{fig:histo_joint}(a,b,c), that the density $\rho =1$ ($\nu =1$) in the gravitational problem is not recovered reliably enough while the value $m=1$ ($\mu =1$) is found quite reliably by the maximum of the histograms Fig.~\ref{fig:histo_joint}(d,e,f). The use of the function \eqref{eq:struct_reg} in the composite residual improves the search for the value $\rho =1$ ($\nu =1$) and somewhat reduces the accuracy of determining the value $m=1$ ($\mu =1$).

\begin{figure}[t]
	\centering
	\includegraphics[width=0.95\linewidth]{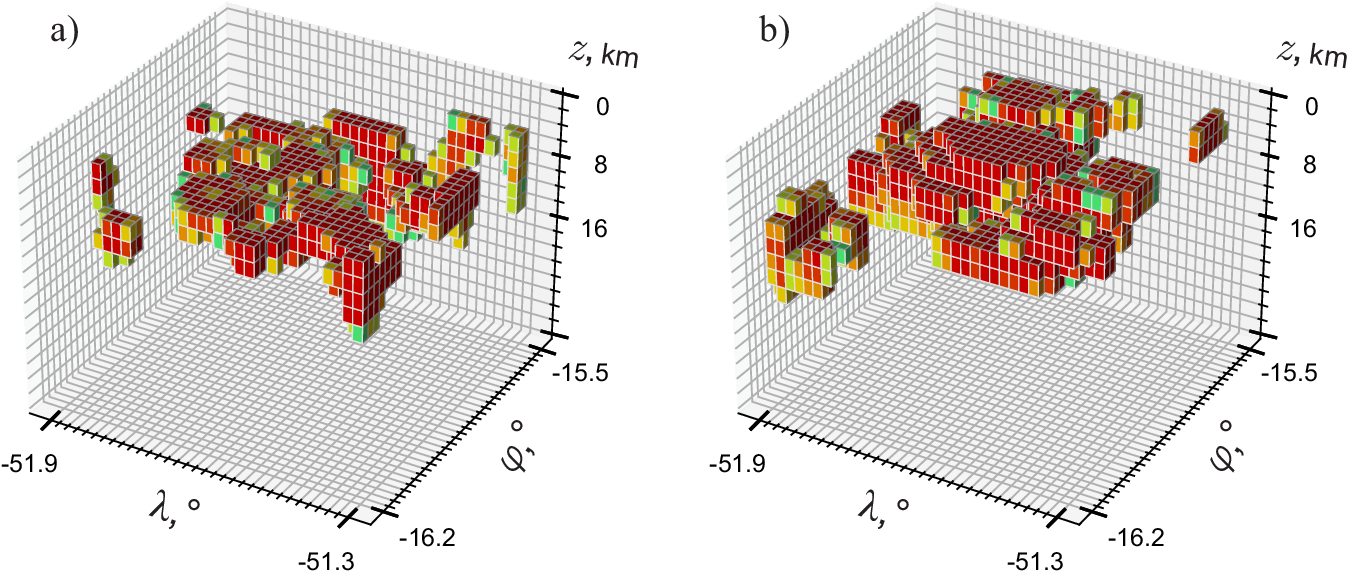}
	\caption{Distribution of gravimetric (left) and magnetic (right) sources during joint processing of experimental data (Jussara region, Goias State, Brazil) for found $\nu$, $\mu$.}
	\label{final}
\end{figure}

This approach can be applied to the joint processing of data from the Jussara region, Goias State, Brazil.
By minimizing function~\eqref{nev_joint}, we obtain, at values $\nu =0.92$, $\mu =1.3$, $\rho =0.92$ g/cm$^3$, $m =1.3$ A/m, the geometric source distributions shown in Fig.~\ref{final}.

\section{Discussion.} \label{sec:discussion}

When considering the results of the work done, let us compare them with related literature sources. We draw special attention to the following important points.

\begin{enumerate}
	\item {\bfseries Dataset comparison.} We used synthetic datasets to train the neural networks, and there is nothing unusual about this.
	Most works employing neural network inversion of geophysical fields use them. For problems of determining the shape of ore bodies, datasets with bodies of simple shape like prisms and steps are often sufficient for training neural networks, which then successfully process real data~\cite{he_cnn21}. As a result of the practical work of such neural networks, body structures are obtained that satisfy expert assessments. As an example, we refer to the well-known NoddyVerse dataset, which implements various geological scenarios. This dataset is the result of many years of practical expert experience in analyzing various geological structures~\cite{noddyverse21}. In our approach, the synthetic dataset is used \emph{for the joint inversion of heterogeneous geophysical fields}.
	
	We note the positive influence of the transfer learning methodology when constructing the dataset and the exceptionally useful introduction of noised data into the dataset. This, in our opinion, is the tool for overcoming (to a certain extent) the instability of approximate solutions.
	
	\item {\bfseries Distinctive features of the neural network for joint inversion.} In a recent work~\cite{bai24}, a typical, thematically similar, approach to the inversion of gravitational fields using neural networks was implemented. Comparing this with our work, one can see a number of significant differences. Firstly, in~\cite{bai24}, only a \emph{two-dimensional} distribution of sources (depending on depth and length of occurrence) was reconstructed from \emph{one-dimensional} observed data. Accordingly, the dataset was built from two-dimensional figures of simpler shape, such as squares and steps. In our approach, three-dimensional bodies of irregular shape, generated randomly from simple bodies, are used. Secondly, there are significant differences in the network architecture. Thirdly, a fundamental difference lies in the fact that in the work~\cite{bai24}, a mean squared error type loss function was used and a so-called PINN-condition was added. We use a loss function of a completely different kind, collectively including several functionals of the Dice type. These functionals separately represent the gravitational and magnetic residuals, and also set the structural residual for the combined assessment of the proximity of gravimetric and magnetic sources.
	
	\item {\bfseries Simultaneous determination of the ore body shape and the quantities $\rho$, $m$ using neural networks in joint inversion.} The conducted numerical experiments on model problems confirmed that for these purposes, composite structural residuals of the form~\eqref{nev_joint} should be used. However, constraints on the quantities $\rho$, $m$ play a significant role in such problems. Inadequate specification of them can lead to significant ambiguity in determining the values of $\rho$, $m$. This can be partially seen in Fig.~\ref{fig:histo_joint}(a,b) for the histograms determining the value of $\rho$.
\end{enumerate}

\section{Conclusions.}

The paper presents a new approach to the joint processing of gravimetric and magnetic data aimed at finding field sources. The approah is based on a two-level neural network approach and the use of special structural residuals with an optimization of neural network elements. The approach provides a significant (by 5 -- 8\%) improvement in the overall residual $\operatorname{Loss}_{\text{joint}}$ between the exact and neural network solutions when processing model data. It can be generalized, both from the point of view of modifying the structural residuals themselves, their combined use and further optimization, and from the point of view of applications to other types of geophysical data. As a possible direction for future work, it is worth noting the generalization of the methodology for using the structural residual to the case of processing full vector magnetic and gravitational data. This also concerns the possibility of simultaneously finding the positions of field sources and the constant values $\rho$, $m$ during joint field inversion.

\textbf{{Acknowledgements.}}
 	The research is carried out using the equipment of the shared research facilities of HPC computing resources at Lomonosov Moscow State University~\cite{Voevodin-2019}.

\textbf{{Funding.}}
  	Russian Science Foundation (RSF-NSFC project 23-41-00002) and National Science Foundation of China NSFC (NSFC-RSF project 12261131494).

\bibliographystyle{plain}
\bibliography{mybibfile}

\end{document}